%% file: multiamen.tex
\documentclass[final,onethmnum, onefignum, onetabnum]{siamltex}
\usepackage{amsmath}
\usepackage{amssymb}
\usepackage[lined,boxruled,linesnumbered]{algorithm2e}
\usepackage{graphicx}

\usepackage{subfigure}

\usepackage{url}

\usepackage{epstopdf}

\usepackage{pgf}
\usepackage{tikz}
\usepackage{color}

\newcommand{\R}{\mathbb{R}}
\newcommand{\calR}{\mathcal{R}}
\newcommand{\calO}{\mathcal{O}}
\newcommand{\calX}{\mathcal{X}}
\newcommand{\calY}{\mathcal{Y}}

\usetikzlibrary{positioning,arrows,automata,fit,shapes}

\title{Multigrid methods combined with low-rank approximation for tensor structured Markov chains} 

\author{Matthias Bolten\footnotemark[1], Karsten Kahl\footnotemark[2], Daniel Kressner\footnotemark[3], Francisco Macedo\footnotemark[3]{\hspace{3pt}}\footnotemark[4] and Sonja Sokolovi\'{c}\footnotemark[2]}

\begin{document}
	\maketitle

	\newcommand{\FF}{F\hspace{-0.075cm}F}
	\footnotetext[1]{Institut f\"ur Mathematik, Universit\"at Kassel, Heinrich-Plett-Str.\ 40, 34132 Kassel, Germany, \texttt{bolten@mathematik.uni-kassel.de}}
	\footnotetext[2]{Fakult\"at f\"ur Mathematik und Naturwissenschaften, Bergische Universit\"at Wuppertal, 42097 Wuppertal, Germany, \texttt{\{kkahl,sokolovic\}@math.uni-wuppertal.de}}
	\footnotetext[3]{EPF Lausanne, SB-MATHICSE-ANCHP, Station 8, CH-1015 Lausanne,
Switzerland, \texttt{\{daniel.kressner,francisco.macedo\}@epfl.ch}}
\footnotetext[4]{IST, Alameda Campus, Av. Rovisco Pais, 1, 1049-001 Lisbon, Portugal}

	\begin{abstract}
Markov chains that describe interacting subsystems suffer, on the one hand, from state space explosion but lead, on the other hand, to highly structured matrices. In this work, we propose a novel tensor-based algorithm to address such tensor structured Markov chains. Our algorithm combines a tensorized multigrid method with AMEn, an optimization-based low-rank tensor solver, for addressing coarse grid problems. Numerical experiments demonstrate that this combination overcomes the limitations incurred when using each of the two methods individually. As a consequence, Markov chain models of unprecedented size from a variety of applications can be addressed.
	\end{abstract}
	
	
	
	\pagestyle{myheadings}
	\thispagestyle{plain}
	\markboth{M. BOLTEN, K. KAHL, D. KRESSNER, F. MACEDO AND S. SOKOLOVI\'{C}}{LOW-RANK TENSOR MULTIGRID FOR MARKOV CHAINS}
	
\section{Introduction}
	
This paper is concerned with the numerical computation of stationary distributions for large-scale continuous--time Markov chains. Mathematically, this task consists of solving the linear system
\begin{equation}\label{eq:Ax0}
Ax=0 \quad \text{with}\quad \mathbf{1}^T x = 1,
\end{equation}
where $A$ is the transposed generator matrix of the Markov chain and $\mathbf{1}$ denotes the vector of all ones. The matrix $A$ is square, nonsymmetric, and satisfies $\mathbf{1}^TA=0$. It is well known~\cite{BermanPlemmons1994} that the irreducibility of $A$ implies existence and uniqueness of the solution of~\eqref{eq:Ax0}.

We specifically consider Markov chains that describe $d$ interacting subsystems. Assuming that the $k$th subsystem has $n_k$ states, the generator matrix usually takes the form
\begin{equation}\label{eq:A}
A = \sum\limits_{t = 1}^T E_1^t \otimes E_2^t \otimes \cdots \otimes E_d^t,
\end{equation}
where $\otimes$ denotes the Kronecker product and $E_k^t \in \R^{n_k \times n_k}$ for $k = 1,\ldots, d$.
Consequently, $A$ has size $n = n_1 n_2 \cdots n_d$, which reflects the fact that 
the states of the Markov chain correspond to all possible combinations of subsystem states. The exponential growth of $n$ with respect to $d$ is
usually called state space explosion~\cite{BuchholzDayar2007}. Applications of models described by~\eqref{eq:A} include queuing theory~\cite{Chan1987,Chan1988,Kaufman1983}, stochastic automata networks~\cite{LangvilleStewart2004, PlateauStewart1997}, analysis of chemical reaction networks~\cite{AndersonCraciunKurtz2010,LevineHwa2007} and telecommunication~\cite{Antunes2005,PhilippeSaadStewart1996}.

The tensor structure of~\eqref{eq:A} can be exploited to yield efficient matrix-vector multiplications in iterative methods for solving~\eqref{eq:Ax0}; see, e.g.,~\cite{LangvilleStewart2004}. However, none of the standard iterative solvers is computationally feasible for larger $d$ because of their need to store vectors of length $n$. To a certain extent, this can be avoided by reducing each $n_k$ with the tensorized multigrid method recently proposed in~\cite{BoltenKahlSokolovic2015}. Still, the need for solving coarse subproblems of size $2^d$ or $3^d$ limits such an approach to modest values of $d$.

Low-rank tensor methods as proposed in~\cite{Buchholz2008b,KressnerMacedo2014} can potentially deal with large values of $d$. The main idea is to view the solution $x$ of~\eqref{eq:A} as an $n_1\times n_2 \times \cdots \times n_d$ tensor and aim at an approximation in a highly compressed, low-rank tensor format. The choice of the format is crucial for the success and practicality of such an approach. In~\cite{Buchholz2008b}, the so called canonical decomposition was used, constituting a natural extension of the concept of product form solutions. Since this format aims at separating all subsystems at the same time, it cannot benefit from an underlying topology and thus often results in relatively large ranks. In contrast, low-rank formats based on tensor networks can be aligned with the topology of interactions between subsystems. In particular, it was demonstrated in~\cite{KressnerMacedo2014} that the so called tensor train  format~\cite{Oseledets2011} appears to be well suited. Alternating optimization techniques are frequently used to obtain approximate solutions within a low-rank tensor format. Specifically,~\cite{KressnerMacedo2014} proposes a variant of the Alternating Minimal Energy method (AMEn)~\cite{Dolgov2014,White2005}. In each step of alternating optimization, a subproblem of the form~\eqref{eq:Ax0} needs to be solved. This turns out to be challenging, although these subproblems are much smaller than the original problem, they are often too large to allow for the solution by a direct method and too ill-conditioned to allow for the solution by an iterative method. It is not known how to design effective preconditioners for such problems.

In this paper, we combine the advantages of two methods. The tensorized multigrid method from~\cite{BoltenKahlSokolovic2015} is used to reduce the mode sizes $n_k$ and the condition number. This, in turn, benefits the use of the low-rank tensor method from~\cite{KressnerMacedo2014} by reducing the size and the condition number of the subproblems.


The rest of this paper is organized as follows. In Section~\ref{sec:tensor} we briefly describe the tensor train format and explain the basic ideas of alternating least squares methods, including AMEn. The tensorized multigrid method is described in Section~\ref{sec:multigrid}. Section~\ref{sec:combination} describes our proposed combination of the tensorized multigrid method with AMEn. In Section~\ref{sec:tests}, the advantages of this combination by a series of numerical experiments involving models from different applications.

\section{Low-rank tensor methods} \label{sec:tensor}

A vector $x \in \R^{n_1 \cdots n_d}$ is turned into a tensor $\calX \in \R^{n_1\times \cdots \times n_d}$ by setting
\begin{equation} \label{eq:tensorindex}
   \calX(i_1, \dots, i_d)=x\big(i_1+(i_2-1)n_1+(i_3-1)n_1n_2+\cdots +(i_d-1)n_1n_2 \cdots n_{d-1}\big)
\end{equation}
with $1 \leq i_{k} \leq n_{k}$ for $k=1,\ldots,d$. In \textsc{Matlab}, this corresponds to the command 
$\texttt{X=reshape(x,n)}$ with $\texttt{n=[n\_1,n\_2,{}\ldots{},n\_d]}$.

\subsection{Tensor train format}\label{sec:tt}

The \emph{tensor train (TT) format} is a multilinear low-dimensional representation of a tensor. Specifically, a tensor $\calX$ is said to be represented in TT format if each entry of the tensor is given by
\begin{equation}\label{eq:tt}
	\calX(i_1,\dots, i_d) = G_1(i_1)\cdot G_2(i_2)\cdots G_d(i_d).
\end{equation} 
The parameter-dependent matrices $G_k(i_k) \in \R^{r_{k-1} \times r_k}$ for $k=1, \ldots, d$
are usually collected in $r_{k-1} \times n_k \times r_k$ tensors, which are called the \emph{TT cores}.
The integers $r_0,r_1,\dots,r_{d-1}, r_d$, with $r_0 = r_d = 1$, determining the sizes of these matrices are called the \emph{TT ranks}.
The complexity of storing $\calX$ in the format~\eqref{eq:tt} is bounded by $(d-2)\widehat{n}\widehat{r}^2+2\widehat{n}\widehat{r}$ if each $n_k \leq \widehat{n}$ and $r_k\leq \widehat{r}$.

For a matrix $A \in \R^{n_1 \cdots n_d \times n_1 \cdots n_d}$, one can define a corresponding \emph{operator TT format} by mapping the row and column indices of $A$ to tensor indices analogous to~\eqref{eq:tensorindex} and letting each entry of $A$ satisfy
\begin{equation}\label{eq:tt_matrix}
 A(i_1,\dots, i_d; j_1,\dots, j_d) = M_1(i_1,j_1)\cdot M_2(i_2,j_2)\cdots M_d(i_d, j_d),
\end{equation}
with parameter-dependent matrices $M_k(i_k,j_k) \in \R^{r_{k-1} \times r_k}$ for $k = 1,\dots, d$. The difference to~\eqref{eq:tt} is that the cores now depend on two parameters instead of one.  A matrix given as a sum of $T$ Kronecker products as in~\eqref{eq:A} can be easily converted into an operator TT format~\eqref{eq:tt_matrix} using, e.g., the techniques described in~\cite{Oseledets2010}. It holds that $r_k \le T$ but often much smaller operator TT ranks can be achieved.

Assuming constant TT ranks, the TT format allows to perform certain elementary operations with a complexity linear (instead of exponential) in $d$. Table~\ref{tab:Costs} summarizes the complexity for operations of interest, which shows that the cost can be expected to be dominated by the TT ranks. For a detailed description of the TT format and its operations, we refer to~\cite{Oseledets2010,Oseledets2011,OseledetsDolgov2012}.
	
	\begin{table}
	\caption{Complexity of operations in TT format for tensors $\calX, \calY \in \R^{n_1 \times \dots \times n_d}$ with TT ranks bounded by $\hat r_\calX$ and $\hat r_\calY$, respectively, and matrix $A \in \R^{(n_1 \times \dots \times n_d) \times (n_1 \times \dots \times n_d)}$ with operator TT ranks bounded by $\hat r_A$. All sizes $n_k$ are bounded by $\hat{n}$.\label{tab:Costs}}
		\begin{center}
			\begin{tabular}{|c|c|c|}\hline
				Operation&Cost& Resulting TT ranks\\ \hline
				Addition of two tensors $\calX+\calY$ & --- & $\hat r_\calX + \hat r_\calY$\\
				Scalar multiplication $\alpha\calX$& $\mathcal{O}(1)$ & $\hat r_\calX$ \\
				Scalar product $\left\langle\calX,\calY\right\rangle$&$\mathcal{O}(d\hat{n}\max\{\hat r_\calX,\hat r_\calY\}^3)$ &---\\
				Matrix-vector product $A\calX$ & $\mathcal{O}(d\hat{n}^2\hat r^2_A\hat r^2_\calX)$& $\hat r_A\hat r_\calX$\\
				Truncation of $\calX$ &$\mathcal{O}(d\widehat{n}\hat r_\calX^3)$& prescribed\\ \hline
			\end{tabular}
				\end{center}
	\end{table}

\subsection{Alternating least squares} \label{sec:als}

In this section, we describe the method of alternating least squares (ALS) from~\cite{KressnerMacedo2014}.

To incorporate the TT format, we first replace~\eqref{eq:Ax0} by the equivalent optimization problem
\begin{equation}\label{eq:min_problem}
	\min \|Ax \| \text{ subject to } \mathbf{1}^Tx=1,
\end{equation}
where $\|\cdot\|$ denotes the Euclidean norm.
We can equivalently view $A$ as a linear operator on $\R^{n_1\times \cdots \times n_d}$ and constrain~\eqref{eq:min_problem} to tensors in TT format:
\begin{equation}\label{eq:const_min_problem}
	\min \|A \calX \| \text{ subject to } \langle \calX, \mathbf{1} \rangle =1,\ \text{$\calX$ is in TT format~\eqref{eq:tt}},
\end{equation}
where $\mathbf{1}$ now refers to the $n_1\times \cdots \times n_d$ tensor of all ones.

Note that the TT format is linear in each of the TT cores. This motivates the use of an alternating least squares (ALS) approach that optimizes the $k$th TT core while keeping all other TT cores fixed. To formulate the subproblem that needs to be solved in each step of ALS, we define the interface matrices
\begin{eqnarray*}
	 	 G_{\leq k-1} &=& \big[G(i_1)\cdots G(i_{k})\big] \in \R^{(n_1\cdots n_{k}) \times r_{k-1}}, \\
	  G_{\geq k+1}&=& \big[G(i_{k +1})\cdots G(i_{d})\big]^T \in \R^{(n_{k+1} \cdots n_d) \times r_{k}}.
\end{eqnarray*}
Without loss of generality, we may assume that the TT format is chosen such that the columns of $G_{\leq k}$
and $G_{\geq k+1}$ are orthonormal; see, e.g.,~\cite{KressnerSteinlechnerUschmajew2013}.
By letting $g_k \in \R^{r_{k-1}n_k r_k}$ contain the vectorization of the $k$th core and setting
\[
	 G_{\neq k} = G_{\leq k-1} \otimes I_{n_{k}} \otimes G_{\geq k+1},
\]
it follows that 
\begin{equation*}
  \text{vec}(\calX) = G_{\neq k} g_k.
\end{equation*}
Inserting this relation into~\eqref{eq:const_min_problem} yields
\[
	\min \|A G_{\neq k} g_k \| \text{ subject to } \langle G_{\neq k} g_k, \mathbf{1} \rangle =1,
\]
which is equivalent to the linear system
\begin{equation}\label{eq:linear_system}
\left[\begin{array}{cc} G_{\neq k}^TA^T\!AG_{\neq k} & \tilde{\mathbf{e}} \\ \tilde{\mathbf{e}}^T & 0 \end{array}\right] \left[\begin{array}{c} g_k \\ \lambda \end{array}\right] = \left[\begin{array}{c} 0 \\ 1 \end{array}\right].
\end{equation}
The vector $\tilde{\mathbf{e}} = G_{\neq k}^T \mathbf e$ can be cheaply computed by tensor contractions. After~\eqref{eq:linear_system} has been solved, the TT format of the tensor $\calX$ is updated by reshaping 
$g_k$ into its $k$th TT core.

One full sweep of ALS consists of applying the described procedure first in a forward sweep over the TT cores $1,2,\ldots,d$ followed
by a backward sweep over the TT cores $d,d-1,\ldots, 1$. 
After each update of a core, an orthogonalization procedure~\cite{Oseledets2011} is applied to ensure the orthonormality of the interface matrices in the subsequent optimization step.

\subsection{AMEn} \label{sec:amen}

The alternating minimal energy (AMEn) method proposed in~\cite{Dolgov2014} for linear systems enriches the TT cores locally by
gradient information, which potentially yields faster convergence than ALS and allows for rank adaptivity. It is sufficient to consider $d = 2$ for illustrating the extension of this procedure to~\eqref{eq:const_min_problem}. The general case $d>2$ then follows analogously to~\cite{Dolgov2014,KressnerSteinlechnerUschmajew2013}
by applying the case $d = 2$ to neighbouring cores.

For $d = 2$, the TT format corresponds to a low-rank factorization $\calX = G_1 G_2^T$ with $G_1 \in \R^{n_1\times r_1}$, $G_2 \in \R^{n_2\times r_2}$. Suppose that the first step of ALS has been performed and $G_1$ has been optimized. We then consider a low-rank approximation of the negative gradient of $\|A \calX\|^2$:
\[
 \calR = -A \calX \approx R_1 R_2^T.
\]
In practice, a rank-2 or rank-3 approximation of $R$ is used. Then the method of steepest descent applied to minimizing $\|A \calX\|^2$ would compute
\[
 \calX + \alpha \calR \approx \begin{pmatrix} G_1 & R_1 \end{pmatrix} \begin{pmatrix} G_2 & \alpha R_2 \end{pmatrix}^T
\]
for some suitably chosen scalar $\alpha$. We now fix (and orthonormalize) the first augmented core $\begin{pmatrix} G_1 & R_1 \end{pmatrix}$. However, instead of using $\begin{pmatrix} G_2 & \alpha R_2 \end{pmatrix}$, we apply the next step of ALS to obtain an optimized second core via the solution of a linear system of the form~\eqref{eq:linear_system}. As a result we obtain an approximation $\calX$ that is at least as good as the one obtained from one forward sweep of ALS without augmentation and, when ignoring the truncation error in $\calR$, at least as good as one step of steepest descent. The described procedure is repeated by augmenting the second core and optimizing the second core, and so on. In each step, the rank of $\calX$ is adjusted by performing low-rank truncation. This rank adaptivity is one of the major advantages of AMEn.

\section{Multigrid}\label{sec:multigrid}

In this section, we recall the multigrid method from~\cite{BoltenKahlSokolovic2015} for solving~\eqref{eq:Ax0} with a matrix $A$ having the tensor structure~\eqref{eq:A}. Special care has to be taken in order to preserve the tensor structure within the multigrid hierarchy. We first introduce the generic components of a multigrid method before explaining the tensor specific construction.

A multigrid approach has the following ingredients: the smoothing scheme, the set of coarse variables, transfer operators (the interpolation operator and the restriction operator) and the coarse grid operator. 

Algorithm~\ref{alg:vcycle} is a prototype of a $V$-cycle and includes the mentioned ingredients. For a detailed description we refer the reader to~\cite{RugeStueben1986,TrottenbergOsterleeSchueller2001}. 

	\begin{algorithm}
		\DontPrintSemicolon
		$ v_\ell = \textnormal{MG}(b_\ell,v_\ell) $\;
		\uIf{ coarsest grid is reached}{solve coarse grid equation $A_\ell v_\ell=b_\ell$.} 
		\Else{Perform $\nu_1$ smoothing steps for $A_\ell v_\ell = b_\ell$ with initial guess $ v_\ell $\; 
			Compute coarse right-hand side  $b_{\ell+1}=Q_\ell(b_\ell-A_\ell v_\ell)$\;
			$ e_{\ell+1}=\textnormal{MG}(b_{\ell+1}, 0)$\;
			 $v_\ell=v_\ell+P_\ell e_{\ell+1}$ \;
			Perform $\nu_2$ smoothing steps for $A_\ell v_\ell = b_\ell$ with initial guess $ v_\ell $\;
		}
		\caption{Multigrid $V$-cycle\label{alg:vcycle}}
	\end{algorithm}

In particular, for a two-grid approach, i.e., $\ell=1,2$, one can describe the realization as follows: the method performs a certain number $\nu_1$ of smoothing steps, using an iterative solver that can be, for instance, weighted Jacobi, Gauss-Seidel or a Krylov subspace method like GMRES~\cite{Saad1996,SaadSchultz1986}; the residual of the current iterate is computed and restricted by a matrix-vector multiplication with the restriction matrix $Q \in \R^{n \times n_c}$; the operator $A_1=A$ is restricted via a Petrov-Galerkin construction to obtain the coarse-grid operator, $A_2=QA_1P \in \R^{n_c\times n_c}$, where $P \in \R^{n_c \times n}$ is the interpolation operator; then we have a recursive call where we solve the coarse grid equation, which is the residual equation; then the error is interpolated and again some smoothing iterations are applied. 

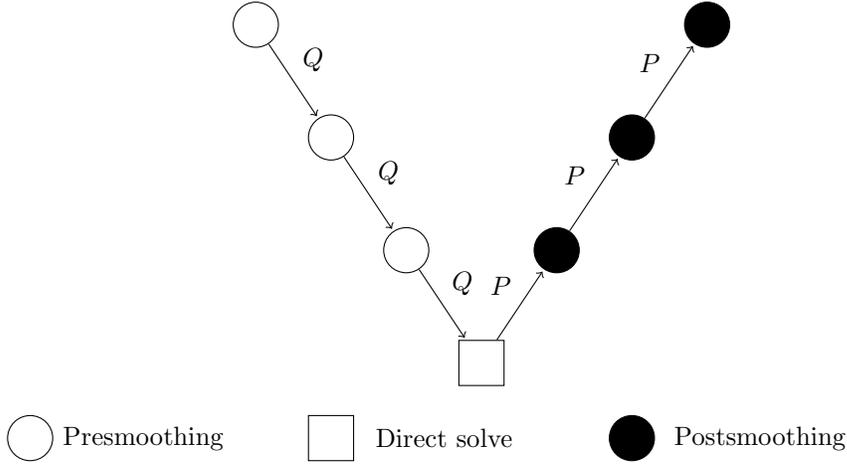
\begin{figure}
		\centering
		\begin{tikzpicture}[shorten >=1pt,auto]
			\draw (2,4.5) node[circle,draw,minimum size = 0.6cm] (A)  {};
			\draw (3,3) node[circle,draw,minimum size = 0.6cm] (B)  {};
			\draw (4,1.5) node[circle,draw,minimum size = 0.6cm] (C)  {};
			\draw (5,0) node[rectangle,draw,minimum size = 0.6cm] (D)  {};
			\draw (6,1.5) node[circle,draw,fill=black,minimum size = 0.6cm] (E) {};
			\draw (7,3) node[circle,draw,fill=black,minimum size = 0.6cm] (F)  {};
			\draw (8,4.5) node[circle,draw,fill=black,minimum size = 0.6cm] (G)  {};
			
			\draw (-1,-1) node[circle,draw,minimum size = 0.6cm] (H)  {};
			\draw (0.5,-1) node (H2)  {Presmoothing}; 
			\draw (3,-1) node[rectangle,draw,minimum size = 0.6cm] (I)  {}; 
			\draw (4.5,-1) node (H2)  {Direct solve};
			\draw (7,-1) node[circle,draw,fill=black,minimum size = 0.6cm] (J)  {};
			\draw (8.7,-1) node (J2)  {Postsmoothing};
			
			\path[->] (A) edge node {$Q$} (B);
			\path[->] (B) edge node {$Q$} (C);
			\path[->] (C) edge node {$Q$} (D);
			\path[->] (D) edge node {$P$} (E);
			\path[->] (E) edge node {$P$} (F);
			\path[->] (F) edge node {$P$} (G);
			
		\end{tikzpicture}
		\caption{Multigrid V-cycle: on each level, a presmoothing iteration is performed before the problem is restricted to the next coarser grid. On the smallest grid, the problem is typically solved exactly by a direct solver. When interpolating back to the finer grids, postsmoothing iterations are applied on each level.}
		\label{fig:vcycle}
	\end{figure}

This $V$-cycle can be performed repeatedly until a certain accuracy of the residual is reached or a maximum number of $V$-cycles have been applied. Instead of stopping at the second grid, because the matrix may still be too large, one can solve the residual equation via a two-grid approach again. By this recursive construction one obtains a multi-level approach, see Fig.~\ref{fig:vcycle}.

No detail has yet been provided on how to choose  $n_c$ and how to obtain the weights for the interpolation and restriction operators $P$ and $Q$. The value $n_c$ is obtained by specifying coarse variables. Geometric coarsening~\cite{TrottenbergOsterleeSchueller2001} or compatible relaxation~\cite{Brandt2000, BrannickFalgout2010} are methods which split the given $n$ variables into fine variables $\mathcal{F}$ and coarse variables $\mathcal{C}$, so that $n=|\mathcal{C}|+|\mathcal{F}|$. If such a splitting is given, $n_c=|C|$, the operators are defined as
	\begin{equation*}
		Q: \R^{|\mathcal{C}\cup\mathcal{F}|} \rightarrow \R^{|\mathcal{C}|},\quad P: \R^{|\mathcal{C}|}\rightarrow\R^{|\mathcal{C}\cup\mathcal{F}|}.
	\end{equation*}
	To obtain the entries for these operators, one can use methods like linear interpolation~\cite{TrottenbergOsterleeSchueller2001} or direct interpolation~\cite{RugeStueben1986,TrottenbergOsterleeSchueller2001}, among others.
	Another approach for choosing a coarse grid is aggregation~\cite{BrezinaManteuffelMcCormichRugeSanders2010}, where one defines a partition of the set of variables and each subset of this partition is associated with one coarse variable. 
	 
In this work we focus on the $V-$cycle strategy. Other strategies, for example $W-$ or $F-$cycles~\cite{TrottenbergOsterleeSchueller2001}, can be applied in a straightforward fashion.
	
\subsection{Tensorized Multigrid}\label{subsec:tensormultigrid}
In order to make Algorithm~\ref{alg:vcycle} applicable to a tensor-structured problem, one has to ensure that the tensor structure is preserved along the multigrid hierarchy. In this, we follow the approach taken in~\cite{BoltenKahlSokolovic2015} and define interpolation and restriction in the following way.

\begin{proposition}\label{pro:QAP}
Let $A$ of the form \eqref{eq:A} be given, with $E_k^t \in \R^{n_k \times n_k}$. Let $P = \bigotimes_{k = 1}^d P_k$ and $Q = \bigotimes_{k = 1}^d Q_k$ with $P_k \in \R^{n_k \times n_k^c}$ and $Q_k \in \R^{n_k^c \times n_k}$ where $n_k^c < n_k$. Then the corresponding Petrov-Galerkin operator satisfies
\[
QAP = \sum\limits_{t = 1}^T\bigotimes_{k = 1}^d Q_k E_k^t P_k.
\]
\end{proposition}

Thus, the task of constructing interpolation and restriction operators becomes a ``local'' task, i.e., each part $P_k$ of the interpolation $P = \bigotimes_{k = 1}^d P_k$ coarsens the $k$th subsystem. In particular, this implies $n_k^{(c)} < n_k$ and the entries of $P_k$ depend largely on the local part of the tensorized operator.

Another important ingredient of the multigrid method is the smoothing scheme. In our setting, it should fulfill two main requirements; it should:
\begin{itemize}
\item[(i)] be applicable to non-symmetric, singular systems;
\item[(ii)] admit an efficient implementation in the TT format. 
\end{itemize}
Requirement (ii) basically means that only the operations listed in Table~\ref{tab:Costs} should be used by the smoother, as most other operations are far more expensive. In this context, one logical choice is GMRES~\cite{Saad1996,SaadSchultz1986} (which also fulfills requirement (i)), which consists of matrix-vector products and orthogonalization steps (i.e., inner products and vector addition). See~\cite{BoltenKahlSokolovic2015} for a discussion of other possible choices for smoothing schemes and their limitations.

\begin{paragraph}{Parameters of the SVD truncation}
We apply the TT-SVD algorithm from~\cite{Oseledets2011} to keep the TT ranks of the iterates in the tensorized multigrid method under control. Except for the application of restriction and interpolation, which both have operator TT rank one by construction, all operations of Algorithm~\ref{alg:vcycle} lead to an increase of the rank of the current iterate.

In particular, truncation has to be performed after line 6 and line 8 of Algorithm~\ref{alg:vcycle}. Concerning the truncation of the restricted residual in line 6, we have observed that we do not need a very strict accuracy to obtain convergence of the global scheme and thus set the value to $10^{-1}$. As for the truncation of the updated iterates $v_\ell$ after line 8, we note that they have highly different norms on the different levels, so that the accuracy for their truncation should depend on the level. Additionally, a dependency on the cycle, following the idea in \cite{KressnerMacedo2014} in which such an adaptive scheme is applied to the sweeps of AMEn, is also included. Precisely, the accuracy depends on the residual norm after the previous cycle. This is motivated by the fact that truncations should be more accurate as we get closer to the desired approximation, while this is not needed while we are still far away from it. Summarizing, the accuracy of the truncation of the different $v_{\ell}$ is thus taken as the norm of $v_{\ell}$ divided by $v_1$ (dependency on the level), times the residual norm after the previous cycle (dependency on the quality of the current approximate solution) times a default value of 10. This ``double'' adaptivity is also used within the GMRES smoother to truncate the occurring vectors.

We also impose a restriction on the maximum TT rank allowed after each truncation. This maximum rank is initially set to $15$ and grows by a factor of $\sqrt{2}$ after each cycle for which the reduction of the residual norm is observed to be smaller than a factor of $\frac{9}{10}$, signalling stagnation. 

\end{paragraph}

\section{Multigrid-AMEn}\label{sec:combination}
 
In Sections~\ref{sec:tensor} and~\ref{sec:multigrid} we have discussed two independent methods for solving~\eqref{eq:Ax0}.
In this section we first discuss the limitations of these two methods and then describe a novel combination that potentially overcomes these limitations.

\subsection{Limitation of AMEn}  \label{sec:limitamen}

Together with orthogonalization and low-rank truncation, one of the computationally most expensive parts of AMEn is the solution of the linear system~\eqref{eq:linear_system}, which has size $r_{k-1} r_k n_k + 1$. A direct solver applied to this linear system has complexity $\calO(\hat r^6 \hat n^3)$
and can thus only be used in the presence of small ranks and mode sizes.

Instead of a direct solver, an iterative solver such as MINRES~\cite{Greenbaum1997,Saad1996} can be applied to~\eqref{eq:linear_system}. The Kronecker structure of $G_{\neq k}^TA^T\!AG_{\neq k}$ inherited by the low operator TT rank of $A$  allows for  efficient matrix-vector multiplications despite the fact that this matrix is not sparse. Unfortunately, we have observed for all the examples considered in Section~\ref{sec:tests} that the condition number of the reduced problem~\eqref{eq:linear_system} grows rapidly as the mode sizes $n_k$ increase. In turn, the convergence of MINRES is severely impaired, often leading to stagnation. It is by no means clear whether it is possible to turn a preconditioner for the original problem into an effective preconditioner for the reduced problem. So far, this has only been achieved via a very particular construction for Laplace-like operators~\cite{KressnerSteinlechnerUschmajew2013}, which is not relevant for the problems under consideration.

\subsection{Limitations of tensorized multigrid} 

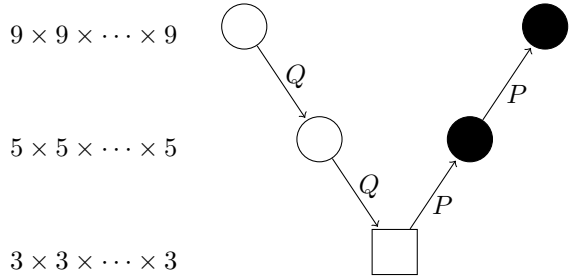
\begin{figure}
\begin{center}
\begin{tikzpicture}[shorten >=1pt,auto]

   \draw (3,3) node[circle,draw,minimum size = 0.6cm] (B)  {};
    \node[label=above:$9 \times 9 \times \dots \times 9$] (1) at (1,2.5) {};
   \draw (4,1.5) node[circle,draw,minimum size = 0.6cm] (C)  {};
      \node[label=above:$5 \times 5 \times \dots \times 5$] (1) at (1,1) {};
   \draw (5,0) node[rectangle,draw,minimum size = 0.6cm] (D)  {};
   \draw (6,1.5) node[circle,draw,fill=black,minimum size = 0.6cm] (E) {};
         \node[label=above:$3 \times 3 \times \dots \times 3$] (1) at (1,-.5) {};
   \draw (7,3) node[circle,draw,fill=black,minimum size = 0.6cm] (F)  {};

   \path[->] (B) edge [above, pos=0.7] node {$\text{ }Q$} (C);
   \path[->] (C) edge [above, pos=0.7] node {$\text{ }Q$} (D);
   \path[->] (D) edge [below, pos=0.6] node {$\text{ }P$} (E);
   \path[->] (E) edge [below, pos=0.6] node {$\text{ }P$} (F);
\end{tikzpicture}
\caption{\label{fig:coarsening}Coarsening process for a problem with mode sizes $9$.}
\end{center}
\end{figure}

The described tensorized multigrid method is limited to modest values of $d$, simply because of the need for solving the problem on the coarsest grid. The size of this problem grows exponentially in $d$. Figure~\ref{fig:coarsening} illustrates the coarsening process if one applies full coarsening to each $E_j^t$ in an overflow queueing problem with mode sizes $9$, as described, e.g., in~\cite[Section 5.1]{BoltenKahlSokolovic2015}; see also Section~\ref{subsec:models} of this paper.  In the case of three levels, a problem of size $3^d$ would need to be addressed by a direct solver on the coarsest grid. Due to the nature of the problem it is not possible to coarse the problem to a single variable in each dimension.

\subsection{Combination of the two methods}

Instead of using a direct method for solving the coarsest-grid system in the tensorized multigrid method, we propose to use AMEn. Due to the fact that the mode sizes on the coarsest grid are small, we expect that it becomes much simpler to solve the reduced problems~\eqref{eq:linear_system} within AMEn.

Note that the problem to be solved on the coarsest grid constitutes a correction equation and thus differs from the original problem~\eqref{eq:Ax0} in having a nonzero right-hand side and incorporating a different linear constraint. To address this problem, we apply AMEn~\cite{Dolgov2014} to the normal equations and ignore the linear constraint. The linear constraint is fixed only at the end of the cycle by explicitly normalizing the obtained approximation, as in \cite{BoltenKahlSokolovic2015}. 

\begin{paragraph}{Parameters of AMEn for the coarsest grid problem} AMEn targets an accuracy that is at the level of the residual from the previous multigrid cycle and we stop AMEn once this accuracy is reached or, at the latest, after $5$ sweeps. A rank-3 approximation of the negative gradient, obtained by ALS as suggested in \cite{Dolgov2014}, is used to augment the cores within AMEn. Reduced problems~\eqref{eq:linear_system} are addressed by a direct solver for size up to $1000$; otherwise  MINRES (without a preconditioner) is used.
\end{paragraph}

\begin{paragraph}{Initial approximation of the solution}
All algorithms are initialized with the tensor that results from solving the coarsest grid problem, using the variant of AMEn described in Section~\ref{sec:amen}, and then bringing it up to the finest level using interpolation, as in \cite{BoltenKahlSokolovic2015}.
\end{paragraph}

\section{Numerical experiments}\label{sec:tests}
In this section, we illustrate the efficiency of our newly proposed algorithm from Section~\ref{sec:combination}. All tests have been performed in \textsc{Matlab} version 2013b, using functions from the \textit{TT-Toolbox}~\cite{Toolbox}.
The execution times have been obtained on a 12-core Intel Xeon CPU X5675, 3.07GHz with 192 GB RAM running 64-Bit Linux version 2.6.32.

\subsection{Model problems}\label{subsec:models}
All benchmark problems used in this paper are taken from the benchmark collection \cite{Macedo2015}, which not only provides a detailed description of the involved matrices but also {\sc Matlab} code.
In total, we consider six different models, which can be grouped into three categories.

\paragraph{Overflow queuing models}

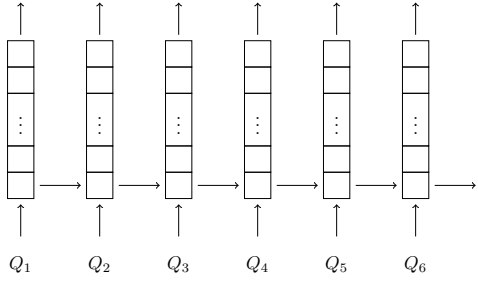
\begin{figure}
\begin{center}
\scalebox{.7}{\input{overflow_tikz.tex}}
\end{center}
\caption{Structure of the model $\mathsf{overflow}$\label{fig:overflow_tikz_new}.}

\end{figure}
The first class of benchmark models consists of the well-known overflow queuing model and two variations thereof. The structure of the model is depicted in Figure~\ref{fig:overflow_tikz_new}. The arrival rates are chosen as $\lambda_k = 1.2 - (k-1)\cdot 0.1$  and the service rates as $\mu_k = 1$ for $k = 1,\dots,d$, as suggested in~\cite{Buchholz2008b}. The variations of the model differ in the interaction between the queues:
\begin{itemize}
\item $\mathsf{overflow}$: Customers which arrive at a full queue try to enter subsequent queues until they find one that is not full. After trying the last queue, they leave the system.
\item $\mathsf{overflowsim}$: As $\mathsf{overflow}$, but customers arriving at a full queue try only one subsequent queue before leaving the system.
\item $\mathsf{overflowpersim}$: As $\mathsf{overflowsim}$, but when the last queue is full, a customer arriving there tries to enter the first queue instead of immediately leaving.
\end{itemize}
For these models, as suggested in \cite{BoltenKahlSokolovic2015}, we choose the interpolation operator $P_k$ as direct interpolation based on the matrices describing the local subsystems, and the restriction operator as its transpose.

\paragraph{Simple tandem queuing network ($\mathsf{kanbanalt2}$)}
\begin{figure}
\begin{center}
\scalebox{.7}{\input{kanbanalt_tikz.tex}}
\caption{Structure of the model $\mathsf{kanbanalt2}$\label{fig:kanbanalt_tikz}.}
\end{center}
\end{figure}
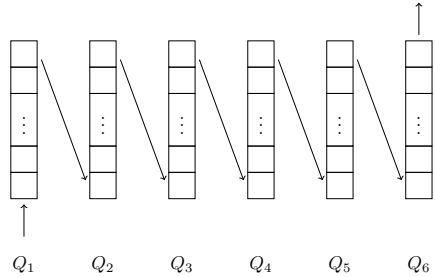
A number $d$ of queues has to be passed through by customers one after the other. Each queue $k$ has its own service rate, denoted by $dep(k)$; and its own capacity, denoted by $cap(k)$. For our tests we choose $dep(k)=1$ for all $k = 1,\dots,d$. The service in queue $k$ can only be finished if queue $k+1$ is not full, so that the served customer can immediately enter the next queue. Customers arrive only at the first queue, with an arrival rate of $1.2$. Figure~\ref{fig:kanbanalt_tikz} illustrates this model.

As only the subsystems corresponding to the first and last dimensions have a non-trivial ``local part'' and the one for the last dimension is associated with a subdiagonal matrix, we construct only $P_1$ via direct interpolation (as in the overflow models) and use linear interpolation for $P_2,\dots,P_d$.

\paragraph{Metabolic pathways}
\begin{figure}
\begin{center}
\subfigure[ ]{\scalebox{.8}{\input{metabolic_tikz.tex}}}
\subfigure[ ]{\scalebox{.8}{\input{metabolic_tikz2.tex}}}
\caption{Structure of the models $\mathsf{directedmetab}$ (a) and $\mathsf{divergingmetab}$ (b)\label{fig:metab_tikz}.}
\end{center}
\end{figure}
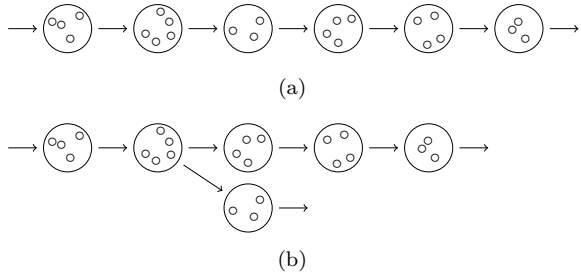
The next model problems we consider come from the field of chemistry, describing stochastic fluctuations in metabolic pathways. In Fig.~\ref{fig:metab_tikz}(a) each node of the given graph describes a metabolite. A flux of substrates can move along the nodes being converted by means of several chemical reactions (an edge between node $k$ and $\ell$ in the graph means that the product of reaction $k$ can be converted further by reaction $\ell$). The rate at which the $k$th reaction happens is given by 
$$\frac{v_k m_k}{m_k + K_k - 1},$$
where $m_k$ is the number of particles of the $k$th substrate and $v_k, K_k$ are constants which we choose as $v_k = 0.1$ and $K_k = 1000$ for all $k = 1,\dots,d$. Note that every substrate $k$ has a maximum capacity of $cap(k)$. This model will be called $\mathsf{directedmetab}$.

$\mathsf{divergingmetab}$ is a variation of this model. Now, one of the metabolites in the reaction network can be converted into two different metabolites, meaning that the reaction path splits into two paths which are independent of each other, as shown in Fig.~\ref{fig:metab_tikz}(b).

The interpolation and restriction operators for these models are chosen in the same way as for $\mathsf{kanbanalt2}$. 

\subsection{Numerical results}\label{subsec:experiments}
In this section, we report the results of the experiments we performed on the models from Section \ref{subsec:models}, in order to compare our proposed method, called ``MultigridAMEn'', to the existing approaches ``AMEn'' and ``Multigrid''.

Throughout all experiments, we stop an iteration when the residual norm $\|Ax\|$ is two orders of magnitude smaller than the residual norm of the tensor of all ones (scaled so that the sum of its entries is one). This happens to be our initial guess for AMEn, but it does not correspond to the initial guesses of Multigrid and MultigridAMEn.

For both multigrid methods, three pre- and postsmoothing steps are applied on each grid.
The number of levels is chosen such that the coarsest grid problem has mode size $3$.

\paragraph{Scaling with respect to the number of subsystems}

\begin{figure}
	\centering
		\subfigure[$\mathsf{overflow}$]{\includegraphics[scale=0.4]{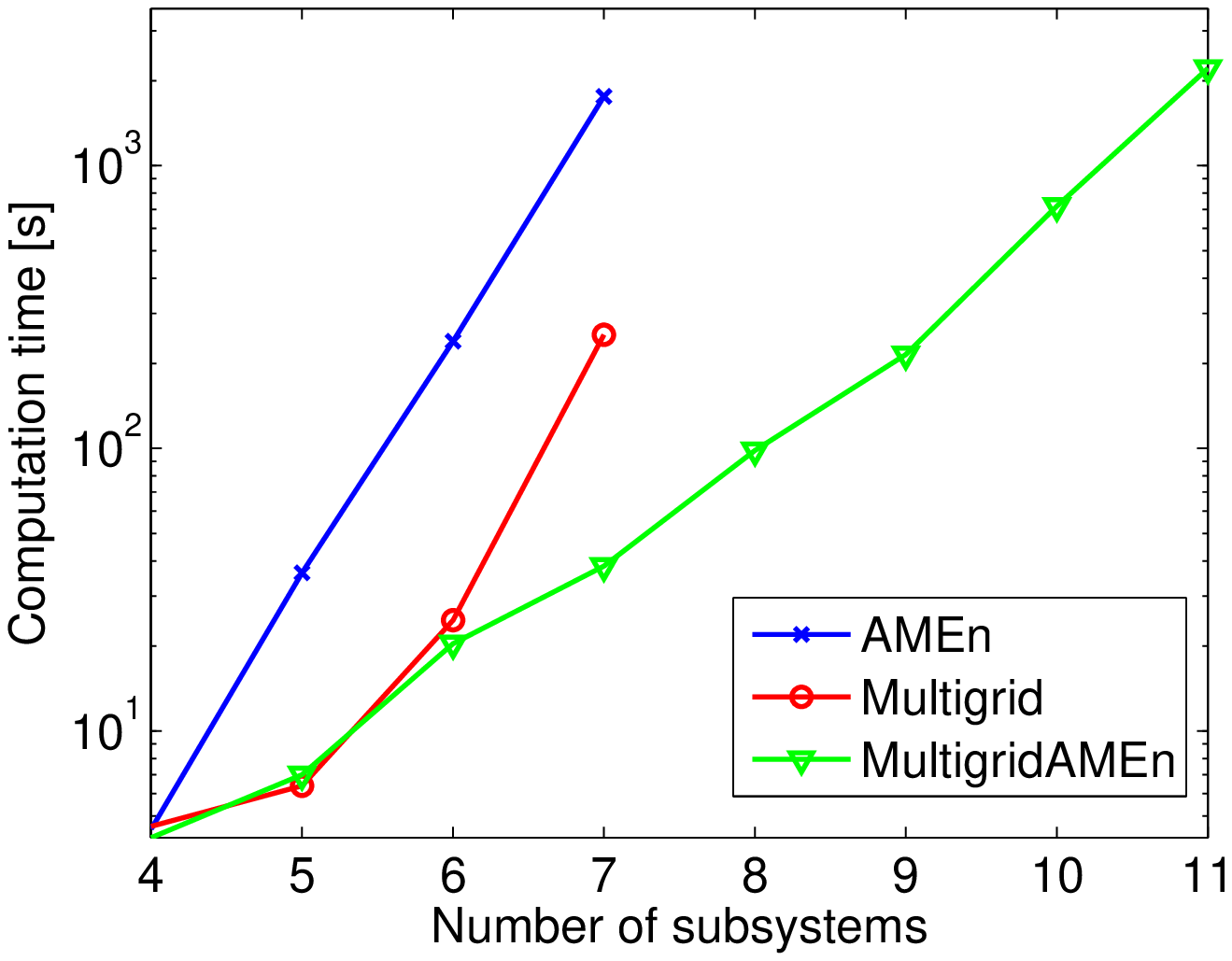}}
		\quad
		\subfigure[$\mathsf{overflowsim}$]{\includegraphics[scale=0.4]{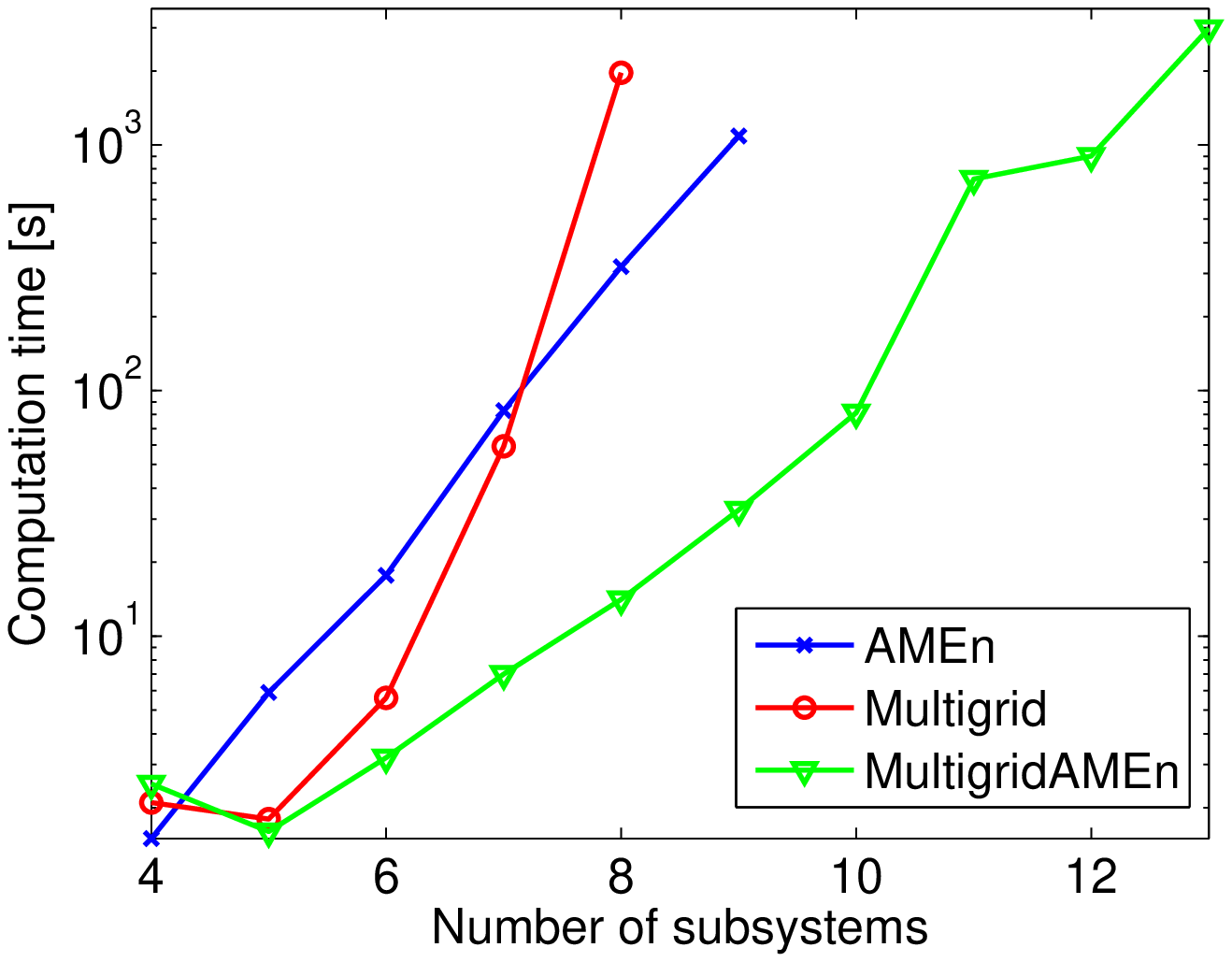}}
		\quad
		\subfigure[$\mathsf{overflowpersim}$]{\includegraphics[scale=0.4]{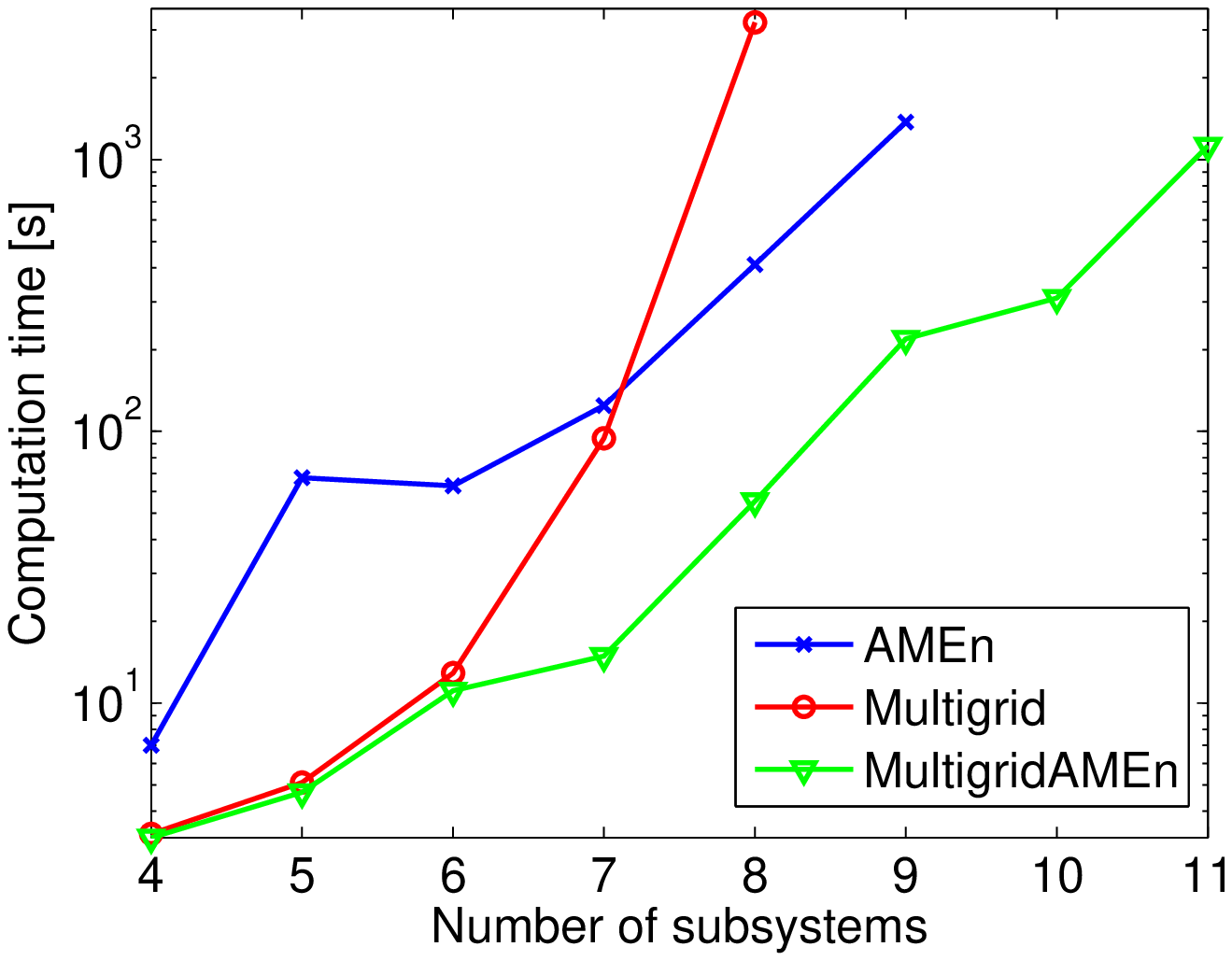}}
		\quad
		\subfigure[$\mathsf{kanbanalt2}$]{\includegraphics[scale=0.4]{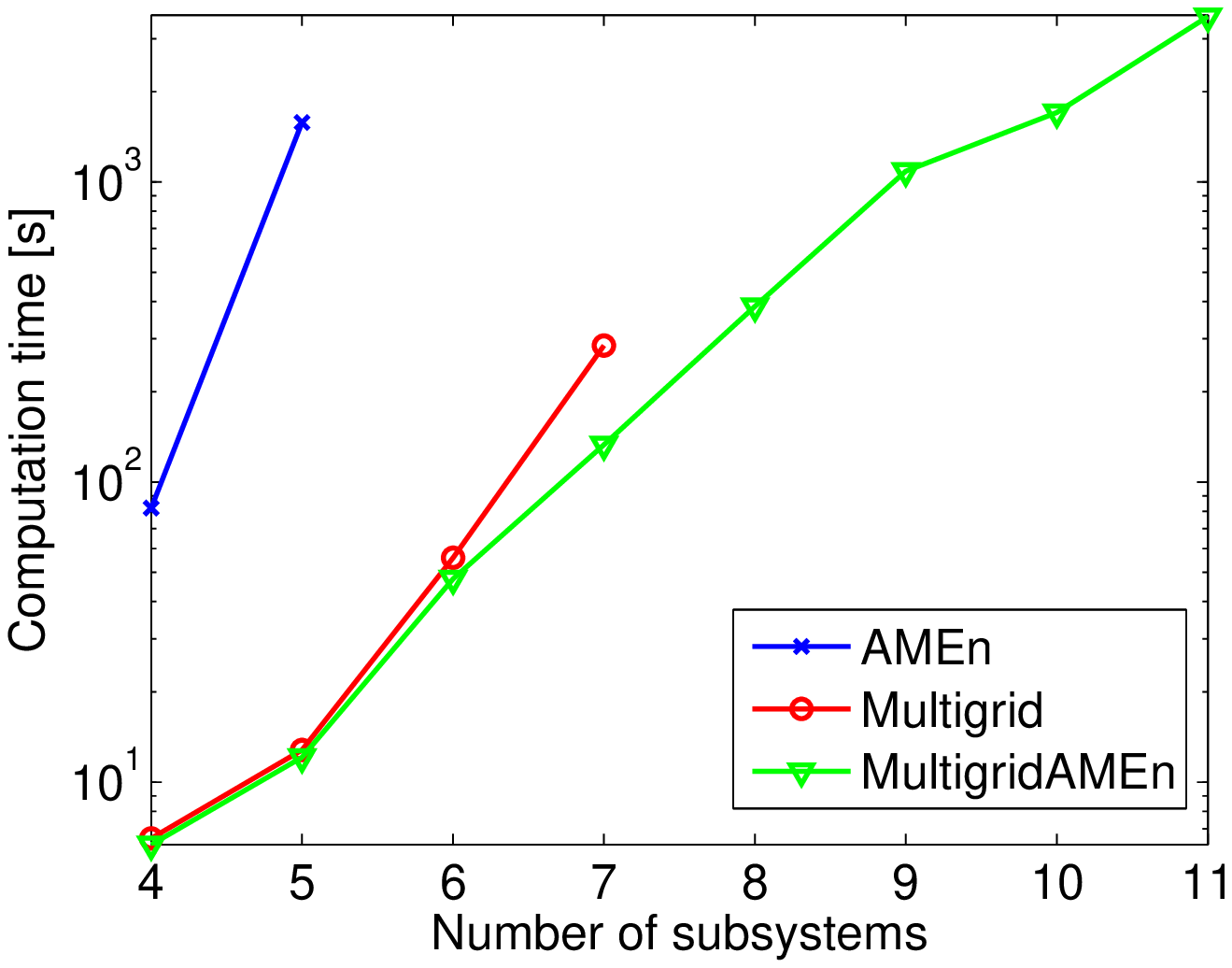}}
		\quad
		\subfigure[$\mathsf{directedmetab}$]{\includegraphics[scale=0.4]{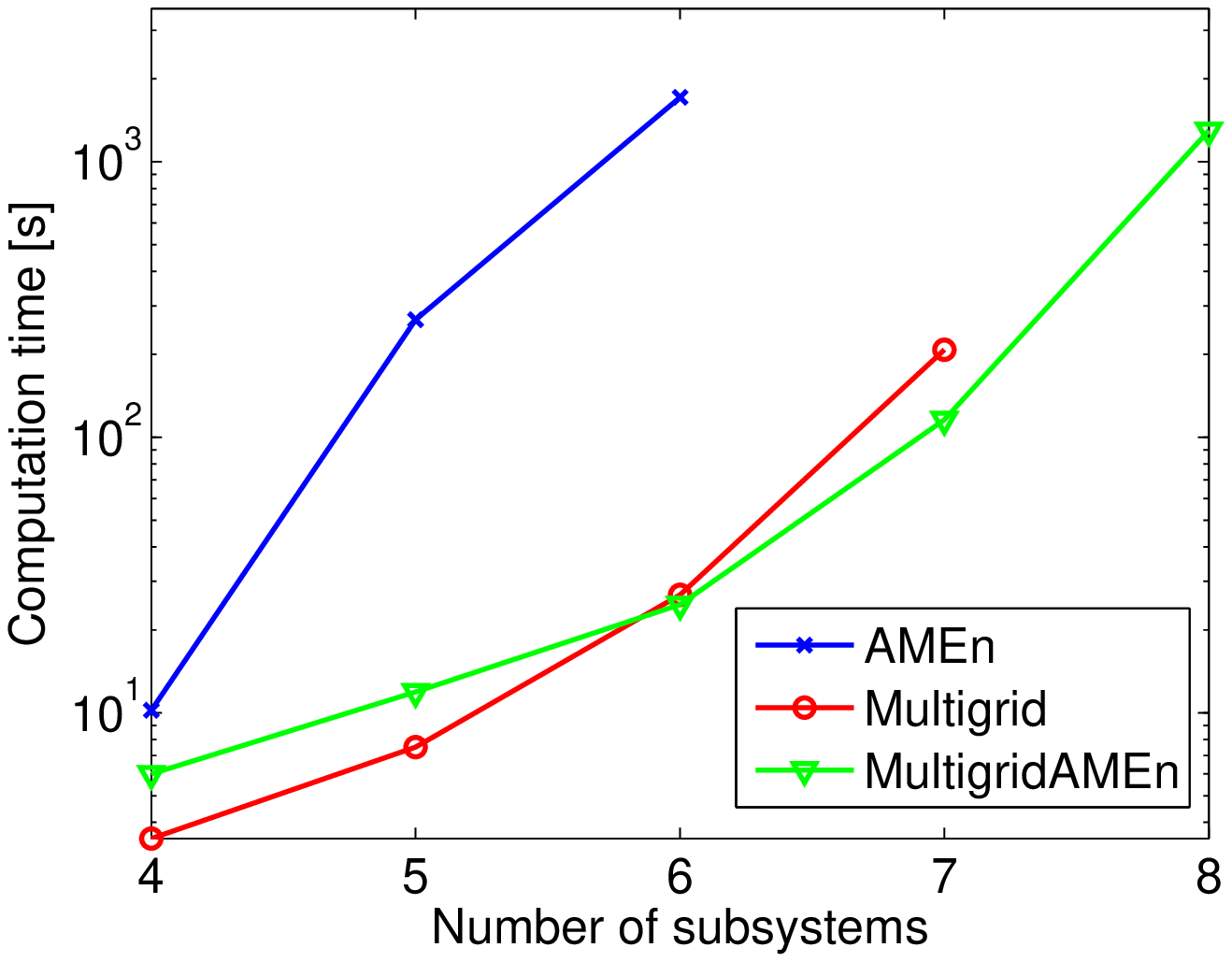}}
		\quad
		\subfigure[$\mathsf{divergingmetab}$]{\includegraphics[scale=0.4]{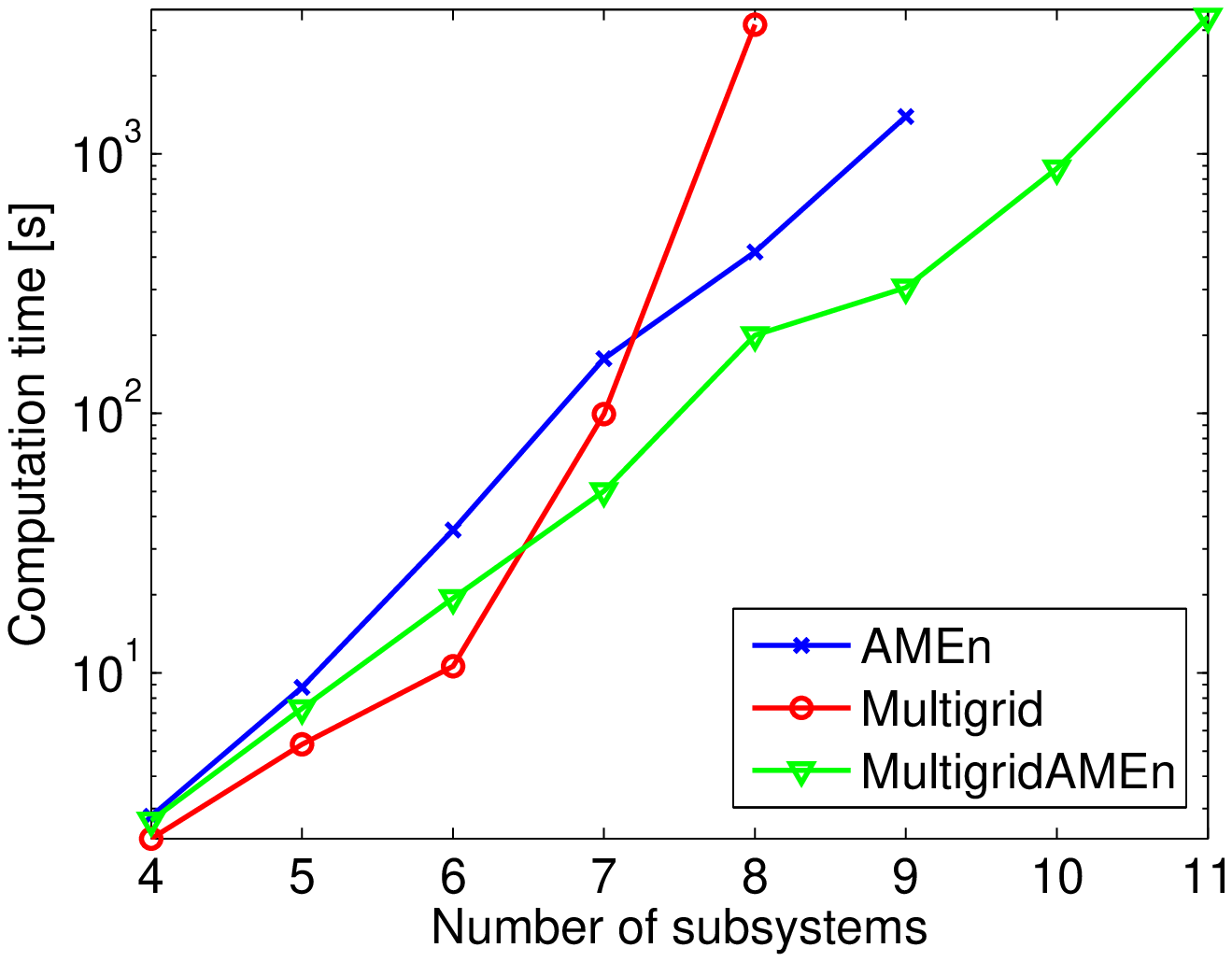}}
		\quad
	\caption{Execution time (in seconds) needed to compute an approximation of the steady state distribution for the benchmark models from Section~\ref{subsec:models}.
All mode sizes are set to $17$.\label{fig:fixedn_new}}
\end{figure}

In order to illustrate the scaling behaviour of the three methods,
we first choose  in all models a capacity of $16$ in each subsystem (i.e., mode sizes 17) and vary $d$, the number of subsystems. Figure~\ref{fig:fixedn_new} displays the obtained execution times.

To provide more insight into the results depicted in Figure~\ref{fig:fixedn_new}, we also give the number of iterations and the maximum rank of the computed approximation for the $\mathsf{overflow}$ model in Table~\ref{tab:fixedn}. For the other models, the observed behaviour is similar and we therefore refrain from providing more detailed data.

\begin{table}
 \caption{Execution time (in seconds), number of iterations, and maximum rank of the computed approximations for $\mathsf{overflow}$ with mode size 17 and varying dimension $d$. The symbol --- indicates that the desired accuracy could not be reached within $3\,600$ seconds.\label{tab:fixedn}}
\footnotesize
\begin{tabular}[c]{c|ccc|ccc|ccc}
&&AMEn&&&Multigrid&&&MultigridAMEn\\
d&time&iter&rank&time&iter&rank&time&iter&rank\\
\hline
4 & 4.5 & 7 & 16 & 4.6& 13 & 13 & 4.2 & 13 & 13 \\
5 & 36.3 & 9 & 23 & 6.4 & 11 & 20 & 7.0 & 11 & 20 \\
6 & 239.4 & 12 & 28 & 24.7 & 17 & 29 & 20.4 & 17 & 29 \\
7 & 1758.4 & 14 & 36 & 252.4 & 24 & 29 & 38.3 & 24 & 29 \\
8 & --- & --- & --- & --- & --- & --- & 98.4 & 28 & 41 \\
9 & --- & --- & --- & --- & --- & --- & 214.8 & 36 & 57 \\
10 & --- & --- & --- & --- & --- & --- & 718.8 & 40 & 80 \\
11 & --- & --- & --- & --- & --- & --- & 2212.2 & 45 & 113 \\
 \end{tabular}
\end{table}

In Figure~\ref{fig:fixedn_new}, we observe that Multigrid and MultigridAMEn behave about the same up to $d = 6$ subsystems. For larger $d$, the cost of solving the coarsest grid problem of size $3^d$ by a direct method becomes prohibitively large within Multigrid. MultigridAMEn is almost always faster than AMEn even for $d = 4$ or $d = 5$. To which extent MultigridAMEn is faster depends on the growth of the TT ranks of the solution with respect to $d$, as these have the largest influence on the performance of AMEn.

Note that the choice of levels in MultigridAMEn is not optimized; it is always chosen such that the coarsest grid mode sizes are three.
We sometimes observed that choosing a larger mode size leads to better performance, but we have not attempted to optimize this choice.

The TT format is a degenerate tree tensor network and thus perfectly matches the topology of interactions in the models {\sf overflowsim}, $\mathsf{kanbanalt2}$, and $\mathsf{directedmetab}$. Compared to $\mathsf{overflowsim}$, the performance is slightly worse for $\mathsf{kanbanalt2}$ and {\sf directedmetab}, possibly because they contain synchronized interactions, that is, interactions associated with a simultaneous change of state in more than one subsystem. In contrast, $\mathsf{overflowsim}$, as well as $\mathsf{overflow}$ and $\mathsf{overflowpersim}$, only have functional interactions, that is, the state of some subsystems determines the rates associated with other subsystems. This seems to be an important factor as the second best performance is observed for $\mathsf{overflowpersim}$, which contains a cycle in the topology of the network and thus does not match the TT format. This robustness with respect to the topology is also reflected by the results for $\mathsf{divergingmetab}$; recall Figure~\ref{fig:metab_tikz}(b).

The maximum problem size that is considered is $17^{13}\approx 9.9 \times 10^{15}$. MultigridAMEn easily deals with larger $d$, but this is the largest configuration for which an execution time below $3\,600$ seconds is obtained. 

\paragraph{Scaling with respect to the mode sizes}

\begin{figure}
	\centering
		\subfigure[$\mathsf{overflow}$]{\includegraphics[scale=0.4]{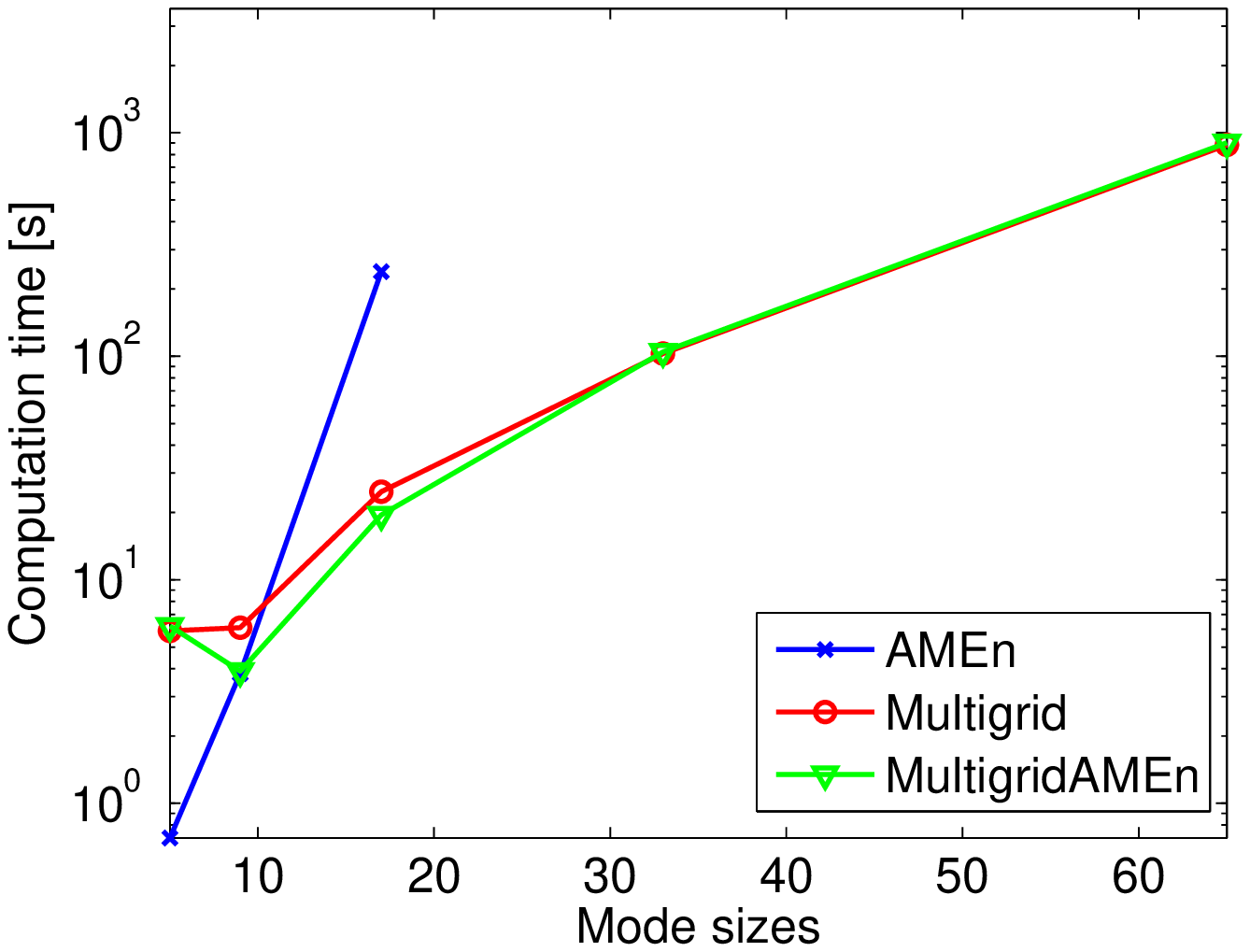}}
		\quad
		\subfigure[$\mathsf{overflowsim}$]{\includegraphics[scale=0.4]{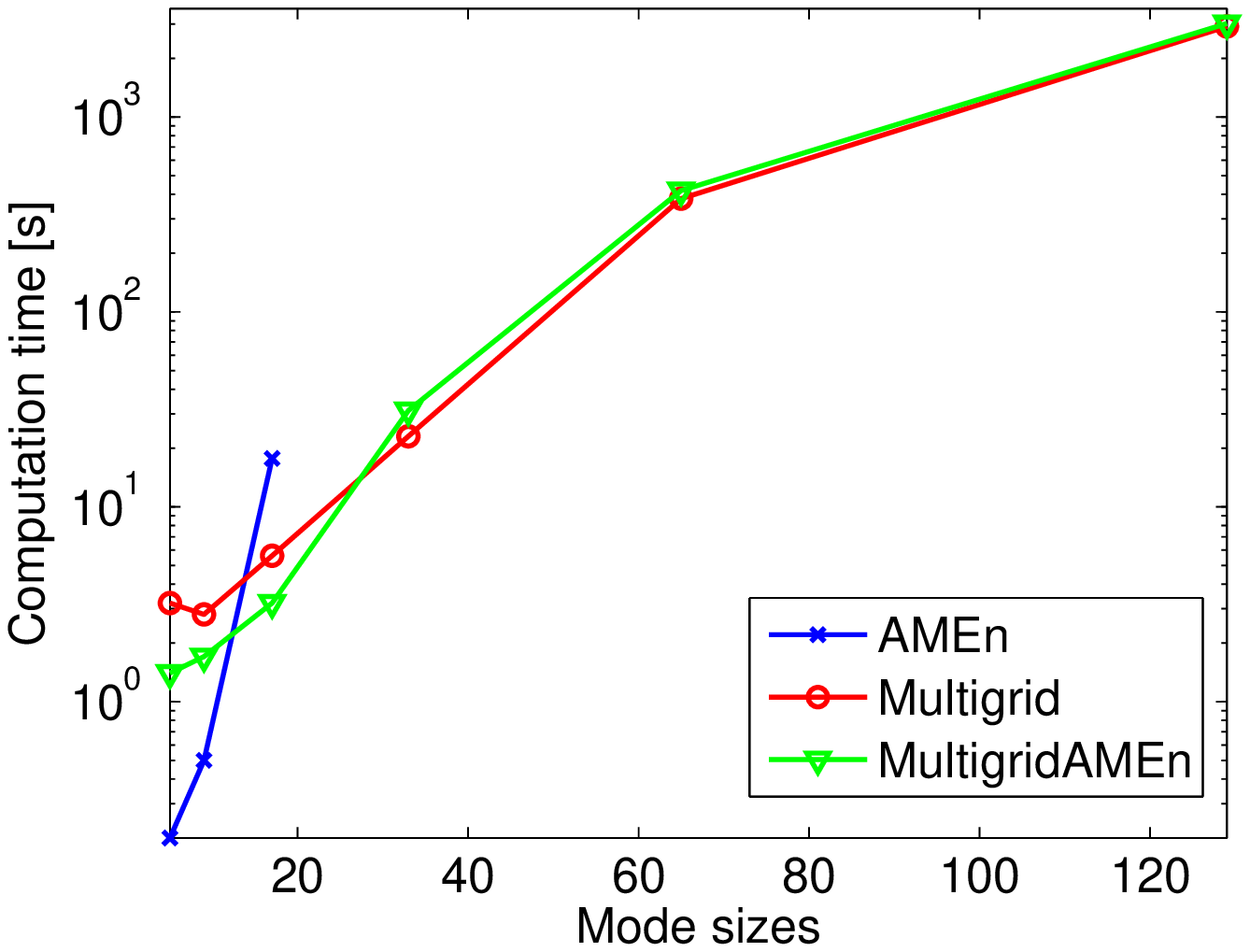}}
		\quad
		\subfigure[$\mathsf{overflowpersim}$]{\includegraphics[scale=0.4]{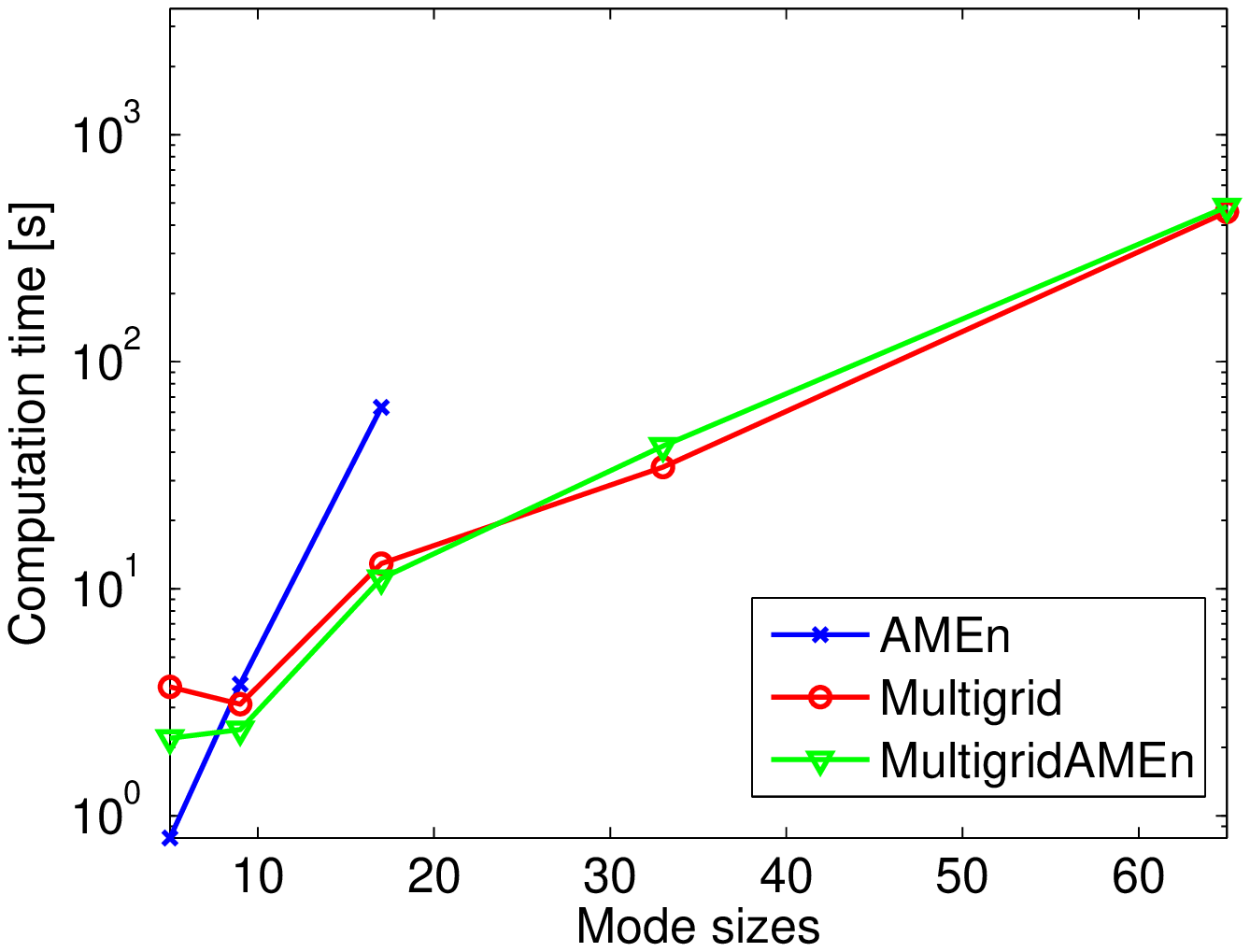}}
		\quad
		\subfigure[$\mathsf{kanbanalt2}$]{\includegraphics[scale=0.4]{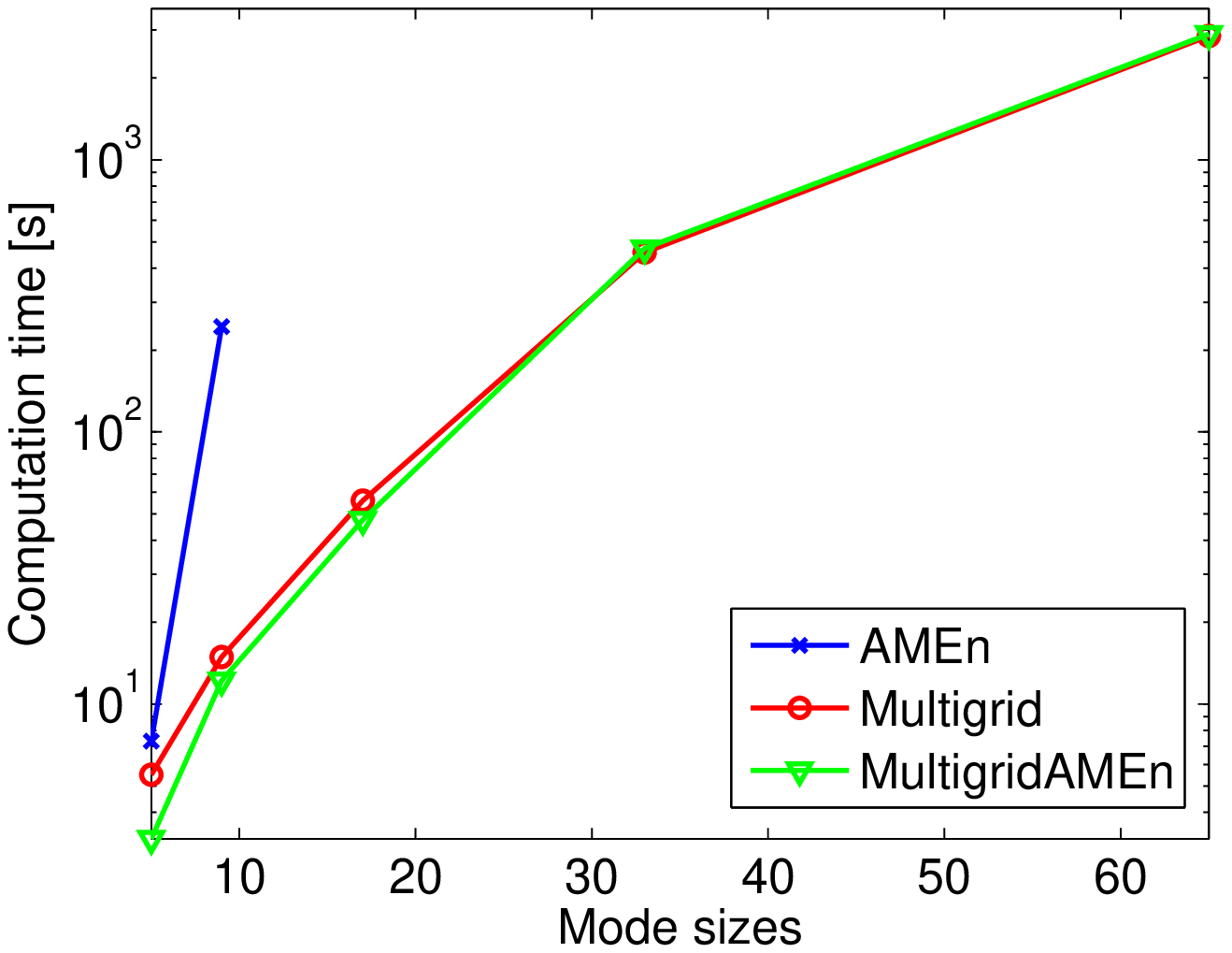}}
		\quad
		\subfigure[$\mathsf{directedmetab}$]{\includegraphics[scale=0.4]{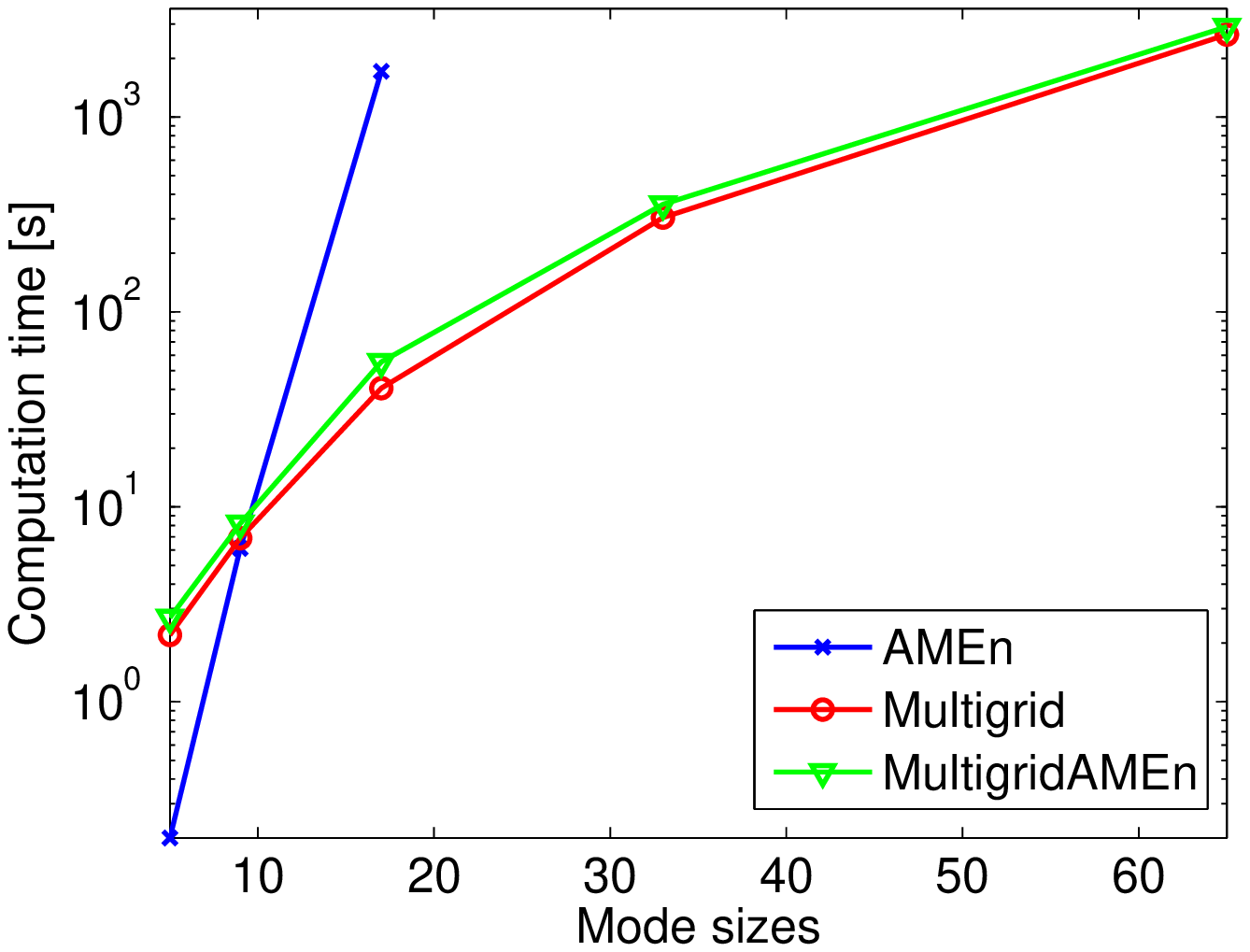}}
		\quad
		\subfigure[$\mathsf{divergingmetab}$]{\includegraphics[scale=0.4]{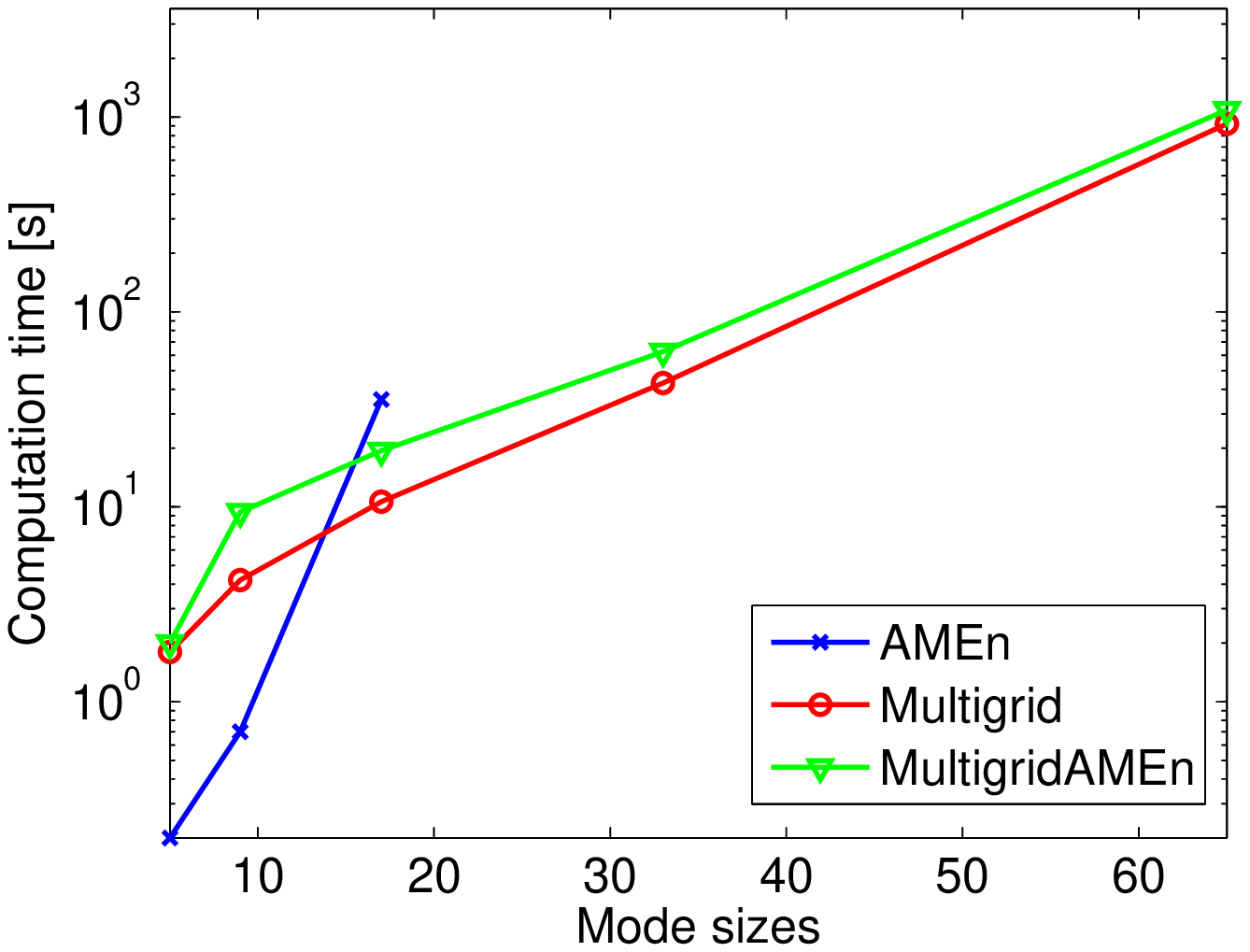}}
		\quad
	\caption{Execution time (in seconds) needed to compute an approximation of the steady state distribution for the benchmark models from Section~\ref{subsec:models}. All models have $d=6$ subsystems.\label{fig:fixedd_new}}
\end{figure}

To also illustrate how the methods scale with respect to increasing mode sizes, we next perform experiments where we fix all models to $d = 6$ subsystems and vary their capacity. The execution times for all models are presented in Figure~\ref{fig:fixedd_new}, while more detailed information for the {\sf overflow} model is given in Table~\ref{tab:fixedd}.

\begin{table}
 \caption{Execution time (in seconds), number of iterations and maximum rank of the computed approximations for $\mathsf{overflow}$ with $d = 6$ and varying mode sizes. The symbol --- indicates that the desired accuracy could not be reached within $3\,600$ seconds.\label{tab:fixedd}}
\small
\begin{tabular}[c]{c|ccc|ccc|ccc}
&&AMEn&&&Multigrid&&&MultiAMEn\\
n&time&iter&rank&time&iter&rank&time&iter&rank\\
\hline
5 & 0.7 & 4 & 13 & 5.9 & 8 & 15 & 6.2 & 8 & 15 \\
9 & 3.8 & 6 & 19 & 6.1 & 8 & 15 & 3.9 & 8 & 15 \\
17 & 239.4& 12 & 28 & 24.8 & 17 & 29 & 19.5 & 17 &29 \\
33 & --- & --- & --- & 102.9 & 17 & 41 & 104.6 & 17 & 41\\
65 & --- & --- & --- & 882.1 & 20 & 57 & 904.1 & 20 & 57\\
 \end{tabular}
\end{table}

Figure~\ref{fig:fixedd_new} shows that AMEn outperforms the two multigrid methods (except for $\mathsf{kanbanalt2}$) for small mode sizes. Depending on the model, the multigrid algorithms start to be faster for mode sizes $9$ or $17$, as the subproblems to be solved in AMEn become too expensive at this point.
The bad performance of AMEn for $\mathsf{kanbanalt2}$ can be explained by the fact that the steady state distribution of this model has rather high TT ranks already for small mode sizes.

Concerning the comparison between the two multigrid methods, no significant difference is visible in Figure~\ref{fig:fixedd_new};
we have already seen in Figure~\ref{fig:fixedn_new} that $d = 6$ is not enough to let the coarsest grid problem solver dominate the computational time in Multigrid.
In fact, Figure~\ref{fig:fixedd_new} nicely confirms that using AMEn for solving the coarsest grid problem does not have an adverse effect on the convergence of multigrid.



The maximum problem size addressed in Figure~\ref{fig:fixedd_new} is $129^6 \approx 4.6 \times 10^{12}$. 

\section{Conclusion}\label{sec:conclusion}

We have proposed a novel combination of two methods, AMEn and Multigrid, for computing the stationary distribution of large-scale tensor structured Markov chains. Our numerical experiments confirm that this combination truly combines the advantages of both methods. As a result, we can address a much wider range of problems in terms of number of subsystems and subsystem states. Also, our experiments demonstrate that the TT format is capable of dealing with a larger variety of applications and topologies compared to what has been previously reported in the literature.

\bibliographystyle{siam}
\bibliography{markov}{}
\end{document}

%% file: overflow_tikz.tex
\begin{tikzpicture}[shorten >=3pt,shorten <= 3pt,auto]

   \draw (0,-1.5) node[] (Q1)  {$Q_1$};
   \draw (0,-1.2) node[] (Q1in)  {};
   \draw (0,3.7) node[] (Q1out)  {};
   \draw (0,0) node[rectangle,draw,minimum width = 0.5cm,minimum height = 0.5cm] (Q11)  {};
   \draw (0,0.5) node[rectangle,draw,minimum width = 0.5cm,minimum height = 0.5cm] (Q12)  {};
   \draw (0,1.25) node[rectangle,draw,minimum width = 0.5cm,minimum height = 1cm] (Q1dots)  {$\vdots$};
   \draw (0,2) node[rectangle,draw,minimum width = 0.5cm,minimum height = 0.5cm] (Q13)  {};
   \draw (0,2.5) node[rectangle,draw,minimum width = 0.5cm,minimum height = 0.5cm] (Q14)  {};
  \path[<-] (Q1out) edge [] node {} (Q14);
   \path[->] (Q1in) edge [] node {} (Q11);

   \draw (1.5,-1.5) node[] (Q2)  {$Q_2$};
   \draw (1.5,-1.2) node[] (Q2in)  {};
   \draw (1.5,3.7) node[] (Q2out)  {};
   \draw (1.5,0) node[rectangle,draw,minimum width = 0.5cm,minimum height = 0.5cm] (Q21)  {};
   \draw (1.5,0.5) node[rectangle,draw,minimum width = 0.5cm,minimum height = 0.5cm] (Q22)  {};
   \draw (1.5,1.25) node[rectangle,draw,minimum width = 0.5cm,minimum height = 1cm] (Q2dots)  {$\vdots$};
   \draw (1.5,2) node[rectangle,draw,minimum width = 0.5cm,minimum height = 0.5cm] (Q23)  {};
   \draw (1.5,2.5) node[rectangle,draw,minimum width = 0.5cm,minimum height = 0.5cm] (Q24)  {};
  \path[<-] (Q2out) edge [] node {} (Q24);
   \path[->] (Q2in) edge [] node {} (Q21);

      \draw (3,-1.5) node[] (Q3)  {$Q_3$};
   \draw (3,-1.2) node[] (Q3in)  {};
   \draw (3,3.7) node[] (Q3out)  {};
   \draw (3,0) node[rectangle,draw,minimum width = 0.5cm,minimum height = 0.5cm] (Q31)  {};
   \draw (3,0.5) node[rectangle,draw,minimum width = 0.5cm,minimum height = 0.5cm] (Q32)  {};
   \draw (3,1.25) node[rectangle,draw,minimum width = 0.5cm,minimum height = 1cm] (Q3dots)  {$\vdots$};
   \draw (3,2) node[rectangle,draw,minimum width = 0.5cm,minimum height = 0.5cm] (Q33)  {};
   \draw (3,2.5) node[rectangle,draw,minimum width = 0.5cm,minimum height = 0.5cm] (Q34)  {};
  \path[<-] (Q3out) edge [] node {} (Q34);
   \path[->] (Q3in) edge [] node {} (Q31);

         \draw (4.5,-1.5) node[] (Q4)  {$Q_4$};
   \draw (4.5,-1.2) node[] (Q4in)  {};
   \draw (4.5,3.7) node[] (Q4out)  {};
   \draw (4.5,0) node[rectangle,draw,minimum width = 0.5cm,minimum height = 0.5cm] (Q41)  {};
   \draw (4.5,0.5) node[rectangle,draw,minimum width = 0.5cm,minimum height = 0.5cm] (Q42)  {};
   \draw (4.5,1.25) node[rectangle,draw,minimum width = 0.5cm,minimum height = 1cm] (Q4dots)  {$\vdots$};
   \draw (4.5,2) node[rectangle,draw,minimum width = 0.5cm,minimum height = 0.5cm] (Q43)  {};
   \draw (4.5,2.5) node[rectangle,draw,minimum width = 0.5cm,minimum height = 0.5cm] (Q44)  {};
  \path[<-] (Q4out) edge [] node {} (Q44);
   \path[->] (Q4in) edge [] node {} (Q41);

         \draw (6,-1.5) node[] (Q5)  {$Q_5$};
   \draw (6,-1.2) node[] (Q5in)  {};
   \draw (6,3.7) node[] (Q5out)  {};
   \draw (6,0) node[rectangle,draw,minimum width = 0.5cm,minimum height = 0.5cm] (Q51)  {};
   \draw (6,0.5) node[rectangle,draw,minimum width = 0.5cm,minimum height = 0.5cm] (Q52)  {};
   \draw (6,1.25) node[rectangle,draw,minimum width = 0.5cm,minimum height = 1cm] (Q5dots)  {$\vdots$};
   \draw (6,2) node[rectangle,draw,minimum width = 0.5cm,minimum height = 0.5cm] (Q53)  {};
   \draw (6,2.5) node[rectangle,draw,minimum width = 0.5cm,minimum height = 0.5cm] (Q54)  {};
  \path[<-] (Q5out) edge [] node {} (Q54);
   \path[->] (Q5in) edge [] node {} (Q51);

         \draw (7.5,-1.5) node[] (Q6)  {$Q_6$};
   \draw (7.5,-1.2) node[] (Q6in)  {};
   \draw (7.5,3.7) node[] (Q6out)  {};
   \draw (7.5,0) node[rectangle,draw,minimum width = 0.5cm,minimum height = 0.5cm] (Q61)  {};
   \draw (7.5,0.5) node[rectangle,draw,minimum width = 0.5cm,minimum height = 0.5cm] (Q62)  {};
   \draw (7.5,1.25) node[rectangle,draw,minimum width = 0.5cm,minimum height = 1cm] (Q6dots)  {$\vdots$};
   \draw (7.5,2) node[rectangle,draw,minimum width = 0.5cm,minimum height = 0.5cm] (Q63)  {};
   \draw (7.5,2.5) node[rectangle,draw,minimum width = 0.5cm,minimum height = 0.5cm] (Q64)  {};
  \path[<-] (Q6out) edge [] node {} (Q64);
   \path[->] (Q6in) edge [] node {} (Q61);   
   
   \draw (9,0) node[rectangle,minimum width = 0.5cm,minimum height = 0.5cm] (Q7)  {};
     
     \path[->] (Q11) edge [] node {} (Q21);  
     \path[->] (Q21) edge [] node {} (Q31);  
     \path[->] (Q31) edge [] node {} (Q41);  
     \path[->] (Q41) edge [] node {} (Q51);
     \path[->] (Q51) edge [] node {} (Q61);  
     \path[->] (Q61) edge [] node {} (Q7);  
   
\end{tikzpicture}

%% file: kanbanalt_tikz.tex
\begin{tikzpicture}[shorten >=3pt,shorten <= 3pt,auto]

   \draw (0,-1.5) node[] (Q1)  {$Q_1$};
   \draw (0,-1.2) node[] (Q1in)  {};
   \draw (0,0) node[rectangle,draw,minimum width = 0.5cm,minimum height = 0.5cm] (Q11)  {};
   \draw (0,0.5) node[rectangle,draw,minimum width = 0.5cm,minimum height = 0.5cm] (Q12)  {};
   \draw (0,1.25) node[rectangle,draw,minimum width = 0.5cm,minimum height = 1cm] (Q1dots)  {$\vdots$};
   \draw (0,2) node[rectangle,draw,minimum width = 0.5cm,minimum height = 0.5cm] (Q13)  {};
   \draw (0,2.5) node[rectangle,draw,minimum width = 0.5cm,minimum height = 0.5cm] (Q14)  {};
   \path[->] (Q1in) edge [] node {} (Q11);
       
   \draw (1.5,-1.5) node[] (Q2)  {$Q_2$};
   \draw (1.5,0) node[rectangle,draw,minimum width = 0.5cm,minimum height = 0.5cm] (Q21)  {};
   \draw (1.5,0.5) node[rectangle,draw,minimum width = 0.5cm,minimum height = 0.5cm] (Q22)  {};
   \draw (1.5,1.25) node[rectangle,draw,minimum width = 0.5cm,minimum height = 1cm] (Q2dots)  {$\vdots$};
   \draw (1.5,2) node[rectangle,draw,minimum width = 0.5cm,minimum height = 0.5cm] (Q23)  {};
   \draw (1.5,2.5) node[rectangle,draw,minimum width = 0.5cm,minimum height = 0.5cm] (Q24)  {};   
   
      \draw (3,-1.5) node[] (Q3)  {$Q_3$};
   \draw (3,0) node[rectangle,draw,minimum width = 0.5cm,minimum height = 0.5cm] (Q31)  {};
   \draw (3,0.5) node[rectangle,draw,minimum width = 0.5cm,minimum height = 0.5cm] (Q32)  {};
   \draw (3,1.25) node[rectangle,draw,minimum width = 0.5cm,minimum height = 1cm] (Q3dots)  {$\vdots$};
   \draw (3,2) node[rectangle,draw,minimum width = 0.5cm,minimum height = 0.5cm] (Q33)  {};
   \draw (3,2.5) node[rectangle,draw,minimum width = 0.5cm,minimum height = 0.5cm] (Q34)  {};   
   
         \draw (4.5,-1.5) node[] (Q4)  {$Q_4$};
   \draw (4.5,0) node[rectangle,draw,minimum width = 0.5cm,minimum height = 0.5cm] (Q41)  {};
   \draw (4.5,0.5) node[rectangle,draw,minimum width = 0.5cm,minimum height = 0.5cm] (Q42)  {};
   \draw (4.5,1.25) node[rectangle,draw,minimum width = 0.5cm,minimum height = 1cm] (Q4dots)  {$\vdots$};
   \draw (4.5,2) node[rectangle,draw,minimum width = 0.5cm,minimum height = 0.5cm] (Q43)  {};
   \draw (4.5,2.5) node[rectangle,draw,minimum width = 0.5cm,minimum height = 0.5cm] (Q44)  {};

         \draw (6,-1.5) node[] (Q5)  {$Q_5$};
   \draw (6,0) node[rectangle,draw,minimum width = 0.5cm,minimum height = 0.5cm] (Q51)  {};
   \draw (6,0.5) node[rectangle,draw,minimum width = 0.5cm,minimum height = 0.5cm] (Q52)  {};
   \draw (6,1.25) node[rectangle,draw,minimum width = 0.5cm,minimum height = 1cm] (Q5dots)  {$\vdots$};
   \draw (6,2) node[rectangle,draw,minimum width = 0.5cm,minimum height = 0.5cm] (Q53)  {};
   \draw (6,2.5) node[rectangle,draw,minimum width = 0.5cm,minimum height = 0.5cm] (Q54)  {};
   
   \draw (7.5,-1.5) node[] (Q6)  {$Q_6$};
   \draw (7.5,3.7) node[] (Q6out)  {};
   \draw (7.5,0) node[rectangle,draw,minimum width = 0.5cm,minimum height = 0.5cm] (Q61)  {};
   \draw (7.5,0.5) node[rectangle,draw,minimum width = 0.5cm,minimum height = 0.5cm] (Q62)  {};
   \draw (7.5,1.25) node[rectangle,draw,minimum width = 0.5cm,minimum height = 1cm] (Q6dots)  {$\vdots$};
   \draw (7.5,2) node[rectangle,draw,minimum width = 0.5cm,minimum height = 0.5cm] (Q63)  {};
   \draw (7.5,2.5) node[rectangle,draw,minimum width = 0.5cm,minimum height = 0.5cm] (Q64)  {};
  \path[<-] (Q6out) edge [] node {} (Q64);   
     
     \draw (0.2,2.75) node[rectangle,minimum width = 0.5cm,minimum height = 0.5cm] (Q1end)  {};
     \draw (1.3,-0.25) node[rectangle,minimum width = 0.5cm,minimum height = 0.5cm] (Q2begin)  {};
     \draw (1.7,2.75) node[rectangle,minimum width = 0.5cm,minimum height = 0.5cm] (Q2end)  {};
     \draw (2.8,-0.25) node[rectangle,minimum width = 0.5cm,minimum height = 0.5cm] (Q3begin)  {};
     \draw (3.2,2.75) node[rectangle,minimum width = 0.5cm,minimum height = 0.5cm] (Q3end)  {};
     \draw (4.3,-0.25) node[rectangle,minimum width = 0.5cm,minimum height = 0.5cm] (Q4begin)  {};
     \draw (4.7,2.75) node[rectangle,minimum width = 0.5cm,minimum height = 0.5cm] (Q4end)  {};
     \draw (5.8,-0.25) node[rectangle,minimum width = 0.5cm,minimum height = 0.5cm] (Q5begin)  {};
     \draw (6.2,2.75) node[rectangle,minimum width = 0.5cm,minimum height = 0.5cm] (Q5end)  {};
     \draw (7.3,-0.25) node[rectangle,minimum width = 0.5cm,minimum height = 0.5cm] (Q6begin)  {};

     \path[->] (Q1end) edge [] node {} (Q2begin);   
     \path[->] (Q2end) edge [] node {} (Q3begin);   
     \path[->] (Q3end) edge [] node {} (Q4begin);   
     \path[->] (Q4end) edge [] node {} (Q5begin);   
     \path[->] (Q5end) edge [] node {} (Q6begin);   
   
\end{tikzpicture}

%% file: metabolic_tikz.tex
\begin{tikzpicture}[shorten >=3pt,shorten <= 3pt,auto]

 \draw (-1.5,0) node[circle,minimum size = 0.8cm] (M0)  {};
   \draw (0,0) node[circle,draw,minimum size = 0.8cm] (M1)  {};
   \draw (0.2,0.2) node[circle,minimum size = 0.05cm] (M11)  {$\circ$};   
   \draw (-0.25,0.1) node[circle,minimum size = 0.05cm] (M12)  {$\circ$};
   \draw (0.05,-0.17) node[circle,minimum size = 0.05cm] (M13)  {$\circ$};   
   \draw (-0.1,0.05) node[circle,minimum size = 0.05cm] (M14)  {$\circ$};   
   
      \draw (1.5,0) node[circle,draw,minimum size = 0.8cm] (M2)  {};
   \draw (1.3,-0.1) node[circle,minimum size = 0.05cm] (M21)  {$\circ$};   
   \draw (1.55,0.27) node[circle,minimum size = 0.05cm] (M22)  {$\circ$};
   \draw (1.72,-0.12) node[circle,minimum size = 0.05cm] (M23)  {$\circ$};   
   \draw (1.69,0.1) node[circle,minimum size = 0.05cm] (M24)  {$\circ$};        
   \draw (1.47,-0.22) node[circle,minimum size = 0.05cm] (M25)  {$\circ$};

        \draw (3,0) node[circle,draw,minimum size = 0.8cm] (M3)  {};
           \draw (2.75,-0.05) node[circle,minimum size = 0.05cm] (M31)  {$\circ$};   
   \draw (3.1,-0.15) node[circle,minimum size = 0.05cm] (M32)  {$\circ$};
   \draw (3.2,0.12) node[circle,minimum size = 0.05cm] (M33)  {$\circ$};

      \draw (4.5,0) node[circle,draw,minimum size = 0.8cm] (M4)  {};
   \draw (4.72,0.15) node[circle,minimum size = 0.05cm] (M41)  {$\circ$};   
   \draw (4.31,-0.1) node[circle,minimum size = 0.05cm] (M42)  {$\circ$};
   \draw (4.48,0.13) node[circle,minimum size = 0.05cm] (M43)  {$\circ$};   
   \draw (4.5,-0.25) node[circle,minimum size = 0.05cm] (M44)  {$\circ$};

         \draw (6,0) node[circle,draw,minimum size = 0.8cm] (M5)  {};
   \draw (6.1,0.21) node[circle,minimum size = 0.05cm] (M51)  {$\circ$};   
   \draw (5.98,-0.27) node[circle,minimum size = 0.05cm] (M52)  {$\circ$};
   \draw (6.2,-0.18) node[circle,minimum size = 0.05cm] (M53)  {$\circ$};   
   \draw (5.82,0.12) node[circle,minimum size = 0.05cm] (M54)  {$\circ$};            
         
      \draw (7.5,0) node[circle,draw,minimum size = 0.8cm] (M6)  {}; 
      \draw (7.61,-0.17) node[circle,minimum size = 0.05cm] (M61)  {$\circ$};   
      \draw (7.39,-0.02) node[circle,minimum size = 0.05cm] (M62)  {$\circ$};
   	  \draw (7.49,0.11) node[circle,minimum size = 0.05cm] (M63)  {$\circ$};  
      
      \draw (9,0) node[circle,minimum size = 0.8cm] (M7)  {};
      
      \path[->] (M0) edge [] node {} (M1); 
     \path[->] (M1) edge [] node {} (M2);   
     \path[->] (M2) edge [] node {} (M3);   
      \path[->] (M3) edge [] node {} (M4);   
     \path[->] (M4) edge [] node {} (M5); 
    \path[->] (M5) edge [] node {} (M6);   
     \path[->] (M6) edge [] node {} (M7);
\end{tikzpicture}

%% file: metabolic_tikz2.tex
\begin{tikzpicture}[shorten >=3pt,shorten <= 3pt,auto]

 \draw (-1.5,0) node[circle,minimum size = 0.8cm] (M0)  {};
   \draw (0,0) node[circle,draw,minimum size = 0.8cm] (M1)  {};
   \draw (0.2,0.2) node[circle,minimum size = 0.05cm] (M11)  {$\circ$};   
   \draw (-0.25,0.1) node[circle,minimum size = 0.05cm] (M12)  {$\circ$};
   \draw (0.05,-0.17) node[circle,minimum size = 0.05cm] (M13)  {$\circ$};   
   \draw (-0.1,0.05) node[circle,minimum size = 0.05cm] (M14)  {$\circ$};   
   
      \draw (1.5,0) node[circle,draw,minimum size = 0.8cm] (M2)  {};
   \draw (1.3,-0.1) node[circle,minimum size = 0.05cm] (M21)  {$\circ$};   
   \draw (1.55,0.27) node[circle,minimum size = 0.05cm] (M22)  {$\circ$};
   \draw (1.72,-0.12) node[circle,minimum size = 0.05cm] (M23)  {$\circ$};   
   \draw (1.69,0.1) node[circle,minimum size = 0.05cm] (M24)  {$\circ$};        
   \draw (1.47,-0.22) node[circle,minimum size = 0.05cm] (M25)  {$\circ$};

        \draw (3,0-1) node[circle,draw,minimum size = 0.8cm] (M3)  {};
           \draw (2.75,-0.05-1) node[circle,minimum size = 0.05cm] (M31)  {$\circ$};   
   \draw (3.1,-0.15-1) node[circle,minimum size = 0.05cm] (M32)  {$\circ$};
   \draw (3.2,0.12-1) node[circle,minimum size = 0.05cm] (M33)  {$\circ$};   
      \draw (4.5,-1) node[circle,minimum size = 0.8cm] (M9)  {};

      \draw (4.5-1.5,0) node[circle,draw,minimum size = 0.8cm] (M4)  {};
   \draw (4.72-1.5,0.15) node[circle,minimum size = 0.05cm] (M41)  {$\circ$};   
   \draw (4.31-1.5,-0.1) node[circle,minimum size = 0.05cm] (M42)  {$\circ$};
   \draw (4.48-1.5,0.13) node[circle,minimum size = 0.05cm] (M43)  {$\circ$};   
   \draw (4.5-1.5,-0.25) node[circle,minimum size = 0.05cm] (M44)  {$\circ$};

         \draw (6-1.5,0) node[circle,draw,minimum size = 0.8cm] (M5)  {};
   \draw (6.1-1.5,0.21) node[circle,minimum size = 0.05cm] (M51)  {$\circ$};   
   \draw (5.98-1.5,-0.27) node[circle,minimum size = 0.05cm] (M52)  {$\circ$};
   \draw (6.2-1.5,-0.18) node[circle,minimum size = 0.05cm] (M53)  {$\circ$};   
   \draw (5.82-1.5,0.12) node[circle,minimum size = 0.05cm] (M54)  {$\circ$};            
         
      \draw (7.5-1.5,0) node[circle,draw,minimum size = 0.8cm] (M6)  {}; 
      \draw (7.61-1.5,-0.17) node[circle,minimum size = 0.05cm] (M61)  {$\circ$};   
      \draw (7.39-1.5,-0.02) node[circle,minimum size = 0.05cm] (M62)  {$\circ$};
   	  \draw (7.49-1.5,0.11) node[circle,minimum size = 0.05cm] (M63)  {$\circ$};  
      
      \draw (9-1.5,0) node[circle,minimum size = 0.8cm] (M7)  {};

      \draw (9,0) node[circle,minimum size = 0.8cm] (M8)  {};
      
      \path[->] (M0) edge [] node {} (M1); 
     \path[->] (M1) edge [] node {} (M2);   
     \path[->] (M2) edge [] node {} (M3);   
      \path[->] (M2) edge [] node {} (M4);   
     \path[->] (M4) edge [] node {} (M5); 
    \path[->] (M5) edge [] node {} (M6);   
     \path[->] (M6) edge [] node {} (M7);
     \path[->] (M3) edge [] node {} (M9); 
\end{tikzpicture}

%% file: multiamen.bbl
\begin{thebibliography}{10}

\bibitem{AndersonCraciunKurtz2010}
{\sc D.~F. Anderson, G.~Craciun, and Th.~G. Kurtz}, {\em Product-form
  stationary distributions for deficiency zero chemical reaction networks},
  Bull.\ Math.\ Biol., 72 (2010), pp.~1947--1970.

\bibitem{Antunes2005}
{\sc N.~Antunes, C.~Fricker, P.~Robert, and D.~Tibi}, {\em Analysis of loss
  networks with routing}, Ann. Appl. Probab., 16 (2006), pp.~2007--2026.

\bibitem{BermanPlemmons1994}
{\sc A.~Berman and R.~J. Plemmons}, {\em Nonnegative Matrices in the
  {M}athematical Sciences}, SIAM, 1994.

\bibitem{BoltenKahlSokolovic2015}
{\sc M.~Bolten, K.~Kahl, and S.~Sokolovi{\'c}}, {\em Multigrid methods for
  tensor structured {M}arkov chains with low rank approximation}, SIAM J. Sci.
  Comput., 38 (2016), pp.~A649--A667.

\bibitem{Brandt2000}
{\sc A.~Brandt}, {\em General highly accurate algebraic coarsening}, Electron.\
  Trans.\ Numer.\ Anal., 10 (2000), pp.~1--20.

\bibitem{BrannickFalgout2010}
{\sc J.~Brannick and R.~Falgout}, {\em Compatible relaxation and coarsening in
  algebraic multigrid}, SIAM J.\ Sci.\ Comput., 32 (2010), pp.~1393--1416.

\bibitem{BrezinaManteuffelMcCormichRugeSanders2010}
{\sc M.~Brezina, T.~A. Manteuffel, S.~F. McCormick, J.~Ruge, and G.~Sanders},
  {\em Towards adaptive smoothed aggregation ($\alpha${SA}) for nonsymmetric
  problems}, SIAM J.\ Sci.\ Comput., 32 (2010), pp.~14--39.

\bibitem{Buchholz2008b}
{\sc P.~Buchholz}, {\em Product form approximations for communicating {M}arkov
  processes}, Perform.\ Eval., 67 (2010), pp.~797--815.

\bibitem{BuchholzDayar2007}
{\sc P.~Buchholz and T.~Dayar}, {\em On the convergence of a class of
  multilevel methods for large sparse {M}arkov chains}, SIAM J.\ Matrix Anal.\
  Appl., 29 (2007), pp.~1025--1049.

\bibitem{Chan1987}
{\sc R.~Chan}, {\em Iterative methods for overflow queueing networks {I}},
  Numer.\ Math., 51 (1987), pp.~143--180.

\bibitem{Chan1988}
\leavevmode\vrule height 2pt depth -1.6pt width 23pt, {\em Iterative methods
  for overflow queueing networks {II}}, Numer.\ Math., 54 (1988), pp.~57--78.

\bibitem{Dolgov2014}
{\sc S.~V. Dolgov and D.~V. Savostyanov}, {\em Alternating minimal energy
  methods for linear systems in higher dimensions}, SIAM J. Sci. Comput., 36
  (2014), pp.~A2248--A2271.

\bibitem{Greenbaum1997}
{\sc A.~Greenbaum}, {\em Iterative Methods for Solving Linear Systems}, Society
  for Industrial and Applied Mathematics, Philadelphia, PA, USA, 1997.

\bibitem{Kaufman1983}
{\sc L.~Kaufman}, {\em Matrix methods for queuing problems}, SIAM J.\ Sci.\
  Statist.\ Comput., 4 (1983), pp.~525--552.

\bibitem{KressnerMacedo2014}
{\sc D.~Kressner and F.~Macedo}, {\em Low-rank tensor methods for communicating
  {M}arkov processes}, in Quantitative Evaluation of Systems, G.~Norman and
  W.~Sanders, eds., vol.~8657 of Lecture Notes in Computer Science, Springer,
  2014, pp.~25--40.

\bibitem{KressnerSteinlechnerUschmajew2013}
{\sc D.~Kressner, M.~Steinlechner, and A.~Uschmajew}, {\em Low-rank tensor
  methods with subspace correction for symmetric eigenvalue problems}, SIAM J.\
  Sci.\ Comput., 36 (2014), pp.~A2346--A2368.

\bibitem{LangvilleStewart2004}
{\sc A.~N. Langville and W.~J. Stewart}, {\em The {K}ronecker product and
  stochastic automata networks}, J.\ Comput.\ Appl.\ Math., 167 (2004),
  pp.~429--447.

\bibitem{LevineHwa2007}
{\sc E.~Levine and T.~Hwa}, {\em Stochastic fluctuations in metabolic
  pathways}, Proc.\ Natl.\ Acad.\ Sci.\ U.S.A., 104 (2007), pp.~9224--9229.

\bibitem{Macedo2015}
{\sc F.~Macedo}, {\em Benchmark problems on stochastic automata networks in
  tensor train format}, tech. report, MATHICSE, EPF Lausanne, Switzerland,
  2015.
\newblock Available from \url{http://anchp.epfl.ch/SAN_TT}.

\bibitem{Oseledets2010}
{\sc I.~V. Oseledets}, {\em Approximation of $2^d\times2^d$ matrices using
  tensor decomposition}, SIAM J.\ Matrix Anal.\ Appl., 31 (2010),
  pp.~2130--2145.

\bibitem{Toolbox}
\leavevmode\vrule height 2pt depth -1.6pt width 23pt, {\em {MATLAB}
  {TT-Toolbox} {V}ersion 2.2}, 2011.
\newblock Available at \url{http://spring.inm.ras.ru/osel/?page\_id=24}.

\bibitem{Oseledets2011}
\leavevmode\vrule height 2pt depth -1.6pt width 23pt, {\em Tensor-{T}rain
  decomposition}, SIAM J.\ Sci.\ Comput., 33 (2011), pp.~2295--2317.

\bibitem{OseledetsDolgov2012}
{\sc I.~V. Oseledets and S.~V. Dolgov}, {\em Solution of linear systems and
  matrix inversion in the {TT}-format}, SIAM J.\ Sci.\ Comput., 34 (2012),
  pp.~A2718--A2739.

\bibitem{PhilippeSaadStewart1996}
{\sc B.~Philippe, Y.~Saad, and W.~J. Stewart}, {\em Numerical methods in
  {M}arkov chain modelling}, Operations Research, 40 (1996), pp.~1156--1179.

\bibitem{PlateauStewart1997}
{\sc B.~Plateau and W.~J. Stewart}, {\em Stochastic automata networks}, in
  Computational Probability, Kluwer Academic Press, 1997, pp.~113--152.

\bibitem{RugeStueben1986}
{\sc J.~Ruge and K.~St\"uben}, {\em Algebraic multigrid}, Multigrid Methods
  (McCormick, S.F., ed.),  (1986).

\bibitem{Saad1996}
{\sc Y.~Saad}, {\em Iterative Methods for Sparse Linear Systems}, Society for
  Industrial and Applied Mathematics, Philadelphia, PA, USA, 2nd~ed., 2003.

\bibitem{SaadSchultz1986}
{\sc Y.~Saad and M.~H. Schultz}, {\em {GMRES}: A generalized minimal residual
  algorithm for solving nonsymmetric linear systems}, SIAM J.\ Sci.\ Comput., 7
  (1986), pp.~856--869.

\bibitem{TrottenbergOsterleeSchueller2001}
{\sc U.~Trottenberg, C.~Osterlee, and A.~Sch\"uller}, {\em Multigrid}, Academic
  Press, 2001.

\bibitem{White2005}
{\sc S.~R. White}, {\em Density matrix renormalization group algorithms with a
  single center site}, Phys.\ Rev.\ B, 72 (2005), p.~180403.

\end{thebibliography}
